\documentclass[12pt]{iopart}
\expandafter\let\csname equation*\endcsname\relax
\expandafter\let\csname endequation*\endcsname\relax
\usepackage[fleqn]{amsmath}
\usepackage{amssymb,amsfonts,amsthm,mathrsfs}
\usepackage{graphicx,multirow,booktabs,textcomp}
\usepackage{float,placeins}
\usepackage{algorithm,algorithmicx,algpseudocode}
\usepackage{enumitem,microtype}
\usepackage{txfonts}
\usepackage[numbers,sort&compress]{natbib}
\usepackage[colorlinks=true,linkcolor=blue,citecolor=blue,urlcolor=blue,hypertexnames=false]{hyperref}
\allowdisplaybreaks
\makeatletter
\newcommand{\namedlabel}[2]{\begingroup\phantomsection\def\@currentlabel{#2}\label{#1}\endgroup}
\makeatother
\theoremstyle{plain}
\newtheorem{theorem}{Theorem}
\newtheorem{proposition}[theorem]{Proposition}
\newtheorem{lemma}{Lemma}

\theoremstyle{definition}

\theoremstyle{remark}

\hypersetup{
  pdftitle={Boundary integral methods and shape reconstruction for periodic fluid-solid interfaces},
  pdfauthor={Yongxin Liu, Junliang Lv, Ning Wang, Yi Wang},
  pdfkeywords={Fluid-solid interaction, Periodic interface, Boundary integral equation, Nystrom method, Shape derivative, Inverse scattering}}

\begin{document}
\title[Boundary integral methods and shape reconstruction]{Boundary integral methods and shape reconstruction for periodic fluid--solid interfaces}
\author{Yongxin Liu$^1$, Junliang Lv$^1$\footnote{Author to whom any correspondence should be addressed.}, Ning Wang$^2$ and Yi Wang$^1$}
\address{$^1$School of Mathematics, Jilin University}
\address{$^2$School of Artificial Intelligence, Hebei Oriental University}
\ead{lvjl@jlu.edu.cn}
\begin{abstract}
This paper studies direct scattering by periodic fluid--solid interfaces and interface reconstruction from multifrequency near-field data. Using Helmholtz decomposition and quasi-periodic single-layer potentials, we formulate a coupled boundary integral system with the standard stress traction, establish its equivalence to the original fluid--solid transmission problem under suitable assumptions, and construct a corresponding Nystr\"om discretization. For near-field observations at finitely many measurement points, we establish local Fr\'echet differentiability of the observation map with a quadratic remainder estimate. In the smooth setting, the derivative transmission problem and its boundary integral formulation are used to compute the observation Jacobian required for interface reconstruction. The interface is reconstructed by regularized Gauss--Newton iteration in a Fourier space that is progressively enriched as the frequency increases. Numerical experiments assess direct-solver accuracy and reconstruction performance for different interfaces, material configurations, and noise levels.
\end{abstract}
\vspace{2pc}
\par
\submitto{\IP}

\section{Introduction}\label{sec1}
An acoustic wave incident on a fluid--solid interface produces a scattered acoustic field in the fluid and excites a transmitted elastic field in the solid. The two fields are coupled through transmission conditions at the interface. Such problems arise in the acoustic characterization of seabed sediments \cite{Kugler2007} and in the excitation and use of ultrasonic surface waves \cite{Declercq2005}. This paper studies direct scattering and near-field reconstruction for unbounded periodic fluid--solid interfaces.

The analysis and computation of direct scattering by periodic fluid--solid interfaces must account for both the acoustic--elastic coupling conditions at the interface and the radiation conditions in the unbounded domains. Hu, Rathsfeld, and Yin use transparent boundary conditions to formulate the problem on a bounded domain within a single periodic cell and establish solvability and numerical convergence results \cite{Yintao2016}. Lin and Lv derive a posteriori estimates that include the finite element error and the acoustic and elastic Dirichlet-to-Neumann (DtN) truncation errors \cite{Linlei2025}. In subsequent joint work, they analyze the well-posedness and exponential convergence of the perfectly matched layer (PML) truncation and derive a posteriori error estimates \cite{LiLinLv2026}. Both studies use these estimates to develop adaptive solution algorithms.

For inverse scattering by periodic fluid--solid interfaces, the far field contains only finitely many propagating Rayleigh modes at a fixed frequency and under fixed illumination, whereas observations close to the interface can retain some evanescent components. Near-field data are therefore often used for interface reconstruction. Sampling methods, including the factorization method, have been developed using acoustic field measurements in the fluid \cite{HuKirschYin2016}. In the three-dimensional biperiodic setting, Cui, Qu, and Wei \cite{CuiQuWei2023} prove, under their stated assumptions, unique determination of the interface from complete scattered acoustic near-field data on a horizontal plane above the interface, corresponding to a countably infinite set of quasi-periodic acoustic incident fields. Domain derivatives and regularization have also been used in iterative methods for reconstructing periodic gratings \cite{Hettlich2002,Hettlich2004}. For random periodic fluid--solid interfaces, Wang, Lin, and Lv \cite{WangLinLv2026} combine a two-step reconstruction method, wavenumber continuation, and Monte Carlo sampling to recover interface statistics. Here we consider deterministic interfaces with known material parameters. Reconstruction uses complex-valued scattered acoustic pressures generated by incident waves with a fixed direction at finitely many frequencies and measured at finitely many points on a horizontal line in the fluid.

Interface reconstruction requires repeated solutions of the direct problem. In homogeneous media, boundary integral methods reduce the field computation to the interface, lowering the dimension by one and avoiding remeshing of the surrounding domains. Integral formulations for scattering by bounded fluid--solid interfaces and their solvability have been studied in \cite{Luke1995,Hsiao2000,Yintao2017}. In two dimensions, Helmholtz decomposition expresses the elastic displacement in terms of two scalar potentials; the two potentials, together with the acoustic field, satisfy a coupled Helmholtz system \cite{Dong2020,SunWangLuChen2023}. Dong et al.\ \cite{Dong2020} use a Nystr\"om discretization and nonlinear integral equations to solve the bounded acoustic--elastic direct problem and reconstruct elastic obstacles from phased or phaseless far-field data. Although their traction formulation differs from the standard stress traction used here, their tangential integration by parts for second-order single-layer potential traces provides a useful reference. Quasi-periodic Green functions incorporate periodicity and radiation conditions into the integral formulation, so that the direct problem can be solved on one period of the interface. The formulation must still accommodate the second-order one-sided traces associated with the standard stress traction, the quasi-periodic densities, and phase-removed kernels suitable for periodic singular quadrature.

The shape derivative describes the linear response of scattering observations to interface changes. Optimization methods based on geometric linearization have been applied to the reconstruction of bounded elastic obstacles \cite{ElschnerHsiaoRathsfeld2009}. For bounded fluid--solid scattering problems, Fr\'echet differentiability with respect to Lipschitz domain deformations has been established \cite{BarucqDerivative2014,Barucq2017}, and a transmission characterization of the derivative has been derived for curvilinear polygonal domains \cite{BarucqMoussaoui2018}. Variational and velocity methods have also been used to study shape derivatives under different boundary regularity assumptions \cite{KangJiang2017}. We differentiate the near-field observation map defined by the coupled direct problem, accounting jointly for the responses of the acoustic and elastic fields to interface perturbations. The observation Jacobian is computed through a linearized transmission problem and a quasi-periodic boundary integral solver.

For periodic fluid--solid scattering with the standard stress traction, we establish a quasi-periodic boundary integral system equivalent to the radiating transmission problem for $C^{3,\sigma}$ interfaces ($0<\sigma<1$) under the non-Wood condition. For smooth interfaces, a kernel splitting that retains the local geometric distance and quasi-periodic phase yields a Nystr\"om discretization with smooth periodic remainders. For near-field observations at finitely many measurement points, we prove local Fr\'echet differentiability of the observation map with respect to $W^{1,\infty}_{\rm per}$ interface perturbations and a quadratic remainder estimate whenever the reference direct problem is uniquely solvable. Under stronger smoothness assumptions, we combine the existing derivative transmission characterization \cite{BarucqMoussaoui2018}, boundary DtN maps, and the integral solver to assemble the shape Jacobian. This Jacobian is used in regularized backtracking Gauss--Newton reconstruction with wavenumber continuation and progressive enrichment of the Fourier space \cite{WangLinLv2026}. Numerical experiments assess direct-solver accuracy and reconstruction performance for different interfaces, material configurations, and noise levels.

Sections~\ref{sec2}--\ref{sec:newton} develop the forward formulation, discretization, shape derivative, and reconstruction method. Section~\ref{sec5} presents numerical experiments, and Section~\ref{sec6} gives conclusions and future work.

\section{Problem formulation}\label{sec2}
\begin{figure}[htbp]\centering
\includegraphics[width=12.38cm]{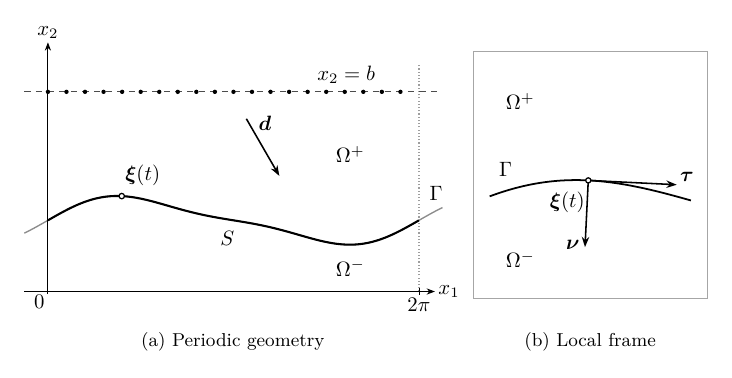}
\caption{A periodic fluid--solid interface and the near-field measurement line. (a) Geometry over one period, with dots marking measurement points on $x_2=b$. (b) Local unit tangent $\boldsymbol\tau$ and unit normal $\boldsymbol\nu$ at $\boldsymbol\xi(t)$, with $\boldsymbol\nu$ pointing into the solid.}\label{fig:model:geometry}
\end{figure}

We consider time-harmonic acoustic scattering by a two-dimensional periodic fluid--solid interface. Let the interface be the graph of a real-valued $2\pi$-periodic function $\gamma\in W^{1,\infty}_{\rm per}(\mathbb R)$, namely,
\[
\Gamma=\{\boldsymbol\xi(t)=(t,\gamma(t))^\top:\ t\in\mathbb R\},
\qquad
S=\{\boldsymbol\xi(t):\ 0\leq t<2\pi\},
\]
where $S$ is one period of the interface and the subscript $\mathrm{per}$ denotes $2\pi$-periodicity. The interface divides the plane into
\[
\Omega^+=\{\boldsymbol x=(x_1,x_2)^\top:\ x_2>\gamma(x_1)\},
\qquad
\Omega^-=\{\boldsymbol x=(x_1,x_2)^\top:\ x_2<\gamma(x_1)\}.
\]
The upper region $\Omega^+$ is occupied by a homogeneous compressible inviscid fluid of density $\rho_a>0$ and sound speed $c>0$, while the lower region $\Omega^-$ is occupied by a homogeneous isotropic elastic solid of density $\rho_e>0$, with Lam\'e constants satisfying $\mu>0$ and $\lambda+\mu>0$. Define the single-period domains $\Omega_0^\pm=\Omega^\pm\cap((0,2\pi)\times\mathbb R)$ and
\[
Q(t)=|\boldsymbol\xi'(t)|=\sqrt{1+[\gamma'(t)]^2},
\qquad
\boldsymbol\nu(t)=\frac{(\gamma'(t),-1)^\top}{Q(t)},
\qquad
\boldsymbol\tau(t)=\frac{(1,\gamma'(t))^\top}{Q(t)}.
\]
The unit normal $\boldsymbol\nu$ points from $\Omega^+$ into $\Omega^-$, and the unit tangent $\boldsymbol\tau=(-\nu_2,\nu_1)^\top$ points in the direction of increasing $t$. The geometry over one period is shown in Fig.~\ref{fig:model:geometry}.

Let $\omega>0$ be the angular frequency. The acoustic, compressional, and shear wavenumbers are, respectively,
\begin{equation}\label{eq:model:wavenumbers}
\kappa_a=\omega/c,\qquad
\kappa_p=\omega\sqrt{\rho_e/(\lambda+2\mu)},\qquad
\kappa_s=\omega\sqrt{\rho_e/\mu}.
\end{equation}
The conditions $\mu>0$ and $\lambda+\mu>0$ imply $0<\kappa_p<\kappa_s$. The time factor $\mathrm e^{-\mathrm i\omega t}$ is suppressed throughout. The plane acoustic wave incident from the fluid side is
\[
p^i(\boldsymbol x)=\mathrm e^{\mathrm i\kappa_a\boldsymbol x\cdot\boldsymbol d},
\qquad
\boldsymbol d=(\cos\theta,-\sin\theta)^\top,
\qquad
\boldsymbol x\in\Omega^+,
\]
where $\theta\in(0,\pi)$ is the angle of incidence measured clockwise from the positive $x_1$-axis, with $\theta=\pi/2$ corresponding to normal incidence. Let $p^s$ and $p=p^i+p^s$ denote the scattered and total acoustic pressures, respectively, and let $\boldsymbol u=(u_1,u_2)^\top$ denote the solid displacement.

The total acoustic pressure $p$ and the displacement $\boldsymbol u$ satisfy the Helmholtz and Navier equations, respectively,
\begin{align}
\Delta p+\kappa_a^2p&=0,\quad \text{in }\Omega^+,\label{Helmholtz}\\
\mu\Delta\boldsymbol u+(\lambda+\mu)\nabla\nabla\!\cdot\boldsymbol u+\omega^2\rho_e\boldsymbol u&=0,\quad \text{in }\Omega^-.\label{Navier}
\end{align}
The strain and stress tensors are defined by
\[
\varepsilon(\boldsymbol w)=\tfrac12\bigl(D\boldsymbol w+(D\boldsymbol w)^\top\bigr),
\qquad
\boldsymbol\sigma(\boldsymbol w)=2\mu\varepsilon(\boldsymbol w)+\lambda(\nabla\!\cdot\boldsymbol w)I,
\]
where $D\boldsymbol w$ is the displacement gradient. The standard stress traction is
\begin{equation}\label{traction}
\mathcal T\boldsymbol w=\boldsymbol\sigma(\boldsymbol w)\boldsymbol\nu
=2\mu\partial_{\boldsymbol\nu}\boldsymbol w+\lambda(\nabla\!\cdot\boldsymbol w)\boldsymbol\nu
-\mu(\operatorname{curl}\boldsymbol w)\boldsymbol\tau,
\end{equation}
where $\partial_{\boldsymbol\nu}\boldsymbol w=(D\boldsymbol w)\boldsymbol\nu$ and $\operatorname{curl}\boldsymbol w=\partial_{x_1}w_2-\partial_{x_2}w_1$. Equation~\eqref{traction} holds for any unit normal.

Traction balance and continuity of the normal velocity across the interface give the transmission conditions
\begin{equation}\label{BDC}
\mathcal T\boldsymbol u=-p\boldsymbol\nu,\qquad
\rho_a\omega^2\boldsymbol u\cdot\boldsymbol\nu=\partial_{\boldsymbol\nu}p,\quad\text{on }\Gamma.
\end{equation}
Set $\alpha=\kappa_a\cos\theta$. The incident field satisfies $p^i(x_1+2\pi,x_2)=\mathrm e^{\mathrm i2\pi\alpha}p^i(x_1,x_2)$, and we impose the corresponding $\alpha$-quasi-periodic conditions
\begin{equation}\label{quasi}
\boldsymbol u(x_1+2\pi,x_2)=\mathrm e^{\mathrm i2\pi\alpha}\boldsymbol u(x_1,x_2),
\qquad
p(x_1+2\pi,x_2)=\mathrm e^{\mathrm i2\pi\alpha}p(x_1,x_2).
\end{equation}
For $s\geq0$, let $H^s_{\alpha,\mathrm{loc}}(\Omega^\pm)$ denote the subspace of $H^s_{\mathrm{loc}}(\Omega^\pm)$ consisting of $\alpha$-quasi-periodic functions, that is, functions for which $\mathrm e^{-\mathrm i\alpha x_1}v$ is $2\pi$-periodic in $x_1$.

Write $\gamma^{\min}=\min_{[0,2\pi)}\gamma$ and $\gamma^{\max}=\max_{[0,2\pi)}\gamma$. To state the radiation conditions and define the quasi-periodic Green function used below, set
\begin{equation}\label{alphazeta}
\alpha_n=\alpha+n,\quad
\zeta_n=\sqrt{\kappa_a^2-\alpha_n^2},\quad
\beta_n=\sqrt{\kappa_p^2-\alpha_n^2},\quad
\eta_n=\sqrt{\kappa_s^2-\alpha_n^2}.
\end{equation}
The square root is taken to be nonnegative on the nonnegative real axis and to have positive imaginary part on the negative real axis. Throughout this paper, we impose the following non-Wood condition.

\noindent\namedlabel{ass:nonwood}{A1}\textbf{(A1)} For all $n\in\mathbb Z$ and $\kappa\in\{\kappa_a,\kappa_p,\kappa_s\}$, $\kappa\neq|\alpha_n|$.

Condition~\eqref{ass:nonwood} excludes grazing modes with zero vertical wavenumber in~\eqref{alphazeta}, so that the Rayleigh representation of the quasi-periodic Green function has no zero denominators.

The scattered acoustic pressure satisfies the upward Rayleigh radiation condition~\cite{Kirsch1993}
\begin{equation}\label{usRayleigh}
p^s(\boldsymbol x)=\sum_{n\in\mathbb Z}p^s_n\mathrm e^{\mathrm i(\alpha_nx_1+\zeta_nx_2)},\qquad x_2>\gamma^{\max}.
\end{equation}
The elastic displacement satisfies the downward Rayleigh radiation condition
\begin{equation}\label{URayleigh}
\boldsymbol u=\sum_{n\in\mathbb Z}\left\{
\mathrm i u_{p,n}\binom{\alpha_n}{-\beta_n}\mathrm e^{\mathrm i(\alpha_nx_1-\beta_nx_2)}
-\mathrm i u_{s,n}\binom{\eta_n}{\alpha_n}\mathrm e^{\mathrm i(\alpha_nx_1-\eta_nx_2)}\right\},\quad x_2<\gamma^{\min},
\end{equation}
where $u_{p,n}$ and $u_{s,n}$ are the modal coefficients of the compressional and shear waves, respectively. In both expansions, only finitely many modes have real vertical wavenumbers; all other modes decay exponentially in their respective radiation directions. For any $a_+>\gamma^{\max}$ and $a_-<\gamma^{\min}$, the right-hand sides of~\eqref{usRayleigh} and~\eqref{URayleigh}, together with all their derivatives, converge uniformly on $x_2\geq a_+$ and $x_2\leq a_-$, respectively.

The governing equations and conditions above define the following direct scattering problem.

\smallskip
\noindent\namedlabel{prob:bvp}{BVP}\textbf{(BVP)} Given $\gamma\in W^{1,\infty}_{\rm per}(\mathbb R)$ and an incident plane wave $p^i$, find the scattered acoustic pressure $p^s\in H^1_{\alpha,\mathrm{loc}}(\Omega^+)$ and displacement $\boldsymbol u\in H^1_{\alpha,\mathrm{loc}}(\Omega^-)^2$, both locally in $H^1$ up to the interface on their respective sides, such that the total acoustic pressure $p=p^i+p^s$ and $\boldsymbol u$ satisfy the governing equations~\eqref{Helmholtz} and~\eqref{Navier}, the quasi-periodic conditions~\eqref{quasi}, the radiation conditions~\eqref{usRayleigh} and~\eqref{URayleigh}, and the transmission conditions~\eqref{BDC}.

Hu, Rathsfeld, and Yin~\cite[Theorem~4.4]{Yintao2016} proved existence of a solution to the corresponding variational problem for a prescribed incident plane wave. Under the material condition $\lambda+2\mu\leq\rho_e c^2$, there exists $\omega_0>0$ such that the solution is unique for $0<\omega\leq\omega_0$; moreover, the problem is uniquely solvable except at a discrete set of frequencies $\mathcal D\subset(0,\infty)$ whose only possible accumulation point is infinity.
Throughout the following analysis, the frequency is chosen so that \textbf{(\ref{prob:bvp})} is uniquely solvable.

\smallskip
The interface regularity assumptions depend on the stage of the analysis.

\noindent\textbf{Regularity conventions.} The weak formulation assumes $\gamma\in W^{1,\infty}_{\rm per}$, and shape differentiability additionally requires perturbations $q\in W^{1,\infty}_{\rm per}$. In Section~\ref{sec3}, the second-order trace formulas and single-layer representation require $\gamma\in C^{2,\sigma}_{\rm per}$, whereas the construction, equivalence, and solvability of the coupled integral system assume $\gamma\in C^{3,\sigma}_{\rm per}$, where $0<\sigma<1$. The kernel decomposition and discretization in Section~\ref{sec4} require $\gamma\in C^\infty_{\rm per}$. The classical transmission characterization of the derivative assumes $\gamma\in C^{3,\sigma}_{\rm per}$ and $q\in C^{2,\sigma}_{\rm per}$, while its boundary integral computation assumes $\gamma,q\in C^\infty_{\rm per}$.

For the near-field inverse problem, the material parameters, incident direction, and observation frequencies $\omega_1<\cdots<\omega_{n_f}$ are assumed known, and the unknown is the $2\pi$-periodic interface function $\gamma$. Fix a measurement height $b$, restrict admissible interfaces to functions satisfying $\gamma^{\max}<b$, and place $M$ equidistant measurement points over one period of the line $x_2=b$. At any fixed frequency, define the near-field observation map, which contains both amplitude and phase information, by
\begin{equation}\label{eq:inv:forward}
\mathcal F(\gamma)=\big(p^s_\gamma(\boldsymbol{x}^{(m)})\big)_{m=0}^{M-1}\in\mathbb C^M,
\quad \boldsymbol{x}^{(m)}=(2\pi m/M,b)^\top.
\end{equation}
The positive distance between the measurement points and the interface, together with interior elliptic regularity, ensures that these point values are well defined. Write the observed data as $\boldsymbol d^\delta=\mathcal F(\gamma_*)+\boldsymbol e$, where $\|\boldsymbol e\|_2\leq\delta$ and $\delta\geq0$ is an upper bound on the noise norm. The near-field inverse problem is to reconstruct the true interface $\gamma_*$ from the noisy complex-valued scattered acoustic pressures measured at the given frequencies.

\section{Boundary integral formulation}\label{sec3}
This section uses the Helmholtz decomposition and quasi-periodic single-layer potentials to derive a boundary integral system for the fluid--solid transmission problem and to establish its equivalence and unique solvability.

\subsection{Quasi-periodic Green function and layer potentials}
The free-space fundamental solution of the two-dimensional Helmholtz equation is $G_0(\boldsymbol x,\boldsymbol y;\kappa)=\allowbreak\mathrm iH_0^{(1)}(\kappa|\boldsymbol x-\boldsymbol y|)/4$, where $H_0^{(1)}$ is the Hankel function of the first kind of order zero. Under the non-Wood condition, the Green function associated with period $2\pi$ and quasi-periodicity parameter $\alpha$ has the Rayleigh representation
\begin{equation}\label{Green2}
G(\boldsymbol{x},\boldsymbol{y};\kappa)
= \frac{\mathrm{i}}{4\pi}
\sum_{n \in \mathbb{Z}}
\frac{1}{\zeta_n(\kappa)}
\mathrm{e}^{\mathrm{i}\{\alpha_n(x_1-y_1)+\zeta_n(\kappa)|x_2-y_2|\}},
\end{equation}
where $\alpha_n=\alpha+n$ and $\zeta_n(\kappa)=\sqrt{\kappa^2-\alpha_n^2}$, with the same square-root branch as in~\eqref{alphazeta}; in particular, $\zeta_n(\kappa_a)$ is the $\zeta_n$ in~\eqref{alphazeta}. We subsequently take $\kappa=\kappa_a,\kappa_p,\kappa_s$.

For any $\boldsymbol z\in S$, translate the periodic cell so that this point lies in its interior. A neighborhood $U_{\boldsymbol z}$ can then be chosen such that
\[
G(\boldsymbol{x},\boldsymbol{y};\kappa)
=
G_0(\boldsymbol{x},\boldsymbol{y};\kappa)
+R_{\boldsymbol z}(\boldsymbol{x},\boldsymbol{y};\kappa),
\qquad
R_{\boldsymbol z}\in C^\infty(U_{\boldsymbol z}\times U_{\boldsymbol z}).
\]
Thus, the jumps across $S$ arise entirely from $G_0$, and the classical jump relations apply directly. For equivalent representations and computation of the Green function, see Linton~\cite{Linton1998}.

We now define the quasi-periodic layer potentials and boundary traces used in the analysis. For $m\in\mathbb N_0$ and $0<\sigma<1$, let $C^{m,\sigma}_\alpha(S)$ denote the space of boundary functions that belong to $C^{m,\sigma}$ after removal of the quasi-periodic factor, that is, $\mathrm e^{-\mathrm i\alpha t}g(\boldsymbol\xi(t))\in C^{m,\sigma}_{\rm per}$. Similarly, $C^\infty_\alpha(S)$ and $H^s_\alpha(S)$ ($s\in\mathbb R$) are defined through pullback and phase removal using the spaces of smooth periodic functions and periodic Sobolev functions, respectively; spaces of negative order are understood by duality.

For a scalar function $v$, write $\mathbf{curl}\,v=(\partial_{x_2}v,-\partial_{x_1}v)^\top$ for its rotated gradient. For a vector $\boldsymbol a=(a_1,a_2)^\top$, write $\boldsymbol a_\perp=(a_2,-a_1)^\top$, so that $\boldsymbol\nu_\perp=-\boldsymbol\tau$.

Define the single-layer potential
\[
V(\boldsymbol{x})=(V_\kappa g)(\boldsymbol x)
=
\int_S G(\boldsymbol{x},\boldsymbol{y};\kappa)
 g(\boldsymbol{y})\,\mathrm d\ell(\boldsymbol y),
\qquad \boldsymbol x\in\mathbb R^2\setminus\Gamma.
\]
In the boundary trace formulas below, integrable singular terms involving $G$ or the normal derivative kernel $\partial_{\boldsymbol\nu(\boldsymbol x)}G$ are integrated in the usual sense. Integrals involving the gradient kernel $\nabla_{\boldsymbol x}G$, the rotated gradient kernel $\mathbf{curl}_{\boldsymbol x}G$, or the tangential derivative kernel $\partial_{\boldsymbol\tau(\boldsymbol x)}G$ are understood in the Cauchy principal value sense; principal value symbols are omitted below. The free-space layer-potential results used here can be found in Chapter~2 of~\cite{ColtonKress1983} and its discussion of the two-dimensional case.

For a field $w$ with boundary traces on both sides of $S$, write
\[
w^\pm(\boldsymbol x)
:=\lim_{\substack{\boldsymbol y\to\boldsymbol x\\ \boldsymbol y\in\Omega^\pm}}
w(\boldsymbol y),
\qquad \boldsymbol x\in S.
\]
The superscripts $+$ and $-$ denote traces from the fluid region $\Omega^+$ and the elastic region $\Omega^-$, respectively.

\begin{proposition}\label{order 1}
Let $0<\sigma<1$, let $\Gamma$ be a $C^{1,\sigma}$ periodic interface, and let $g\in C^{0,\sigma}_\alpha(S)$. Then $V$ is continuous across $S$, and
\[
(\nabla V)^\pm(\boldsymbol{x})
=
\int_S
\nabla_{\boldsymbol{x}}G(\boldsymbol{x},\boldsymbol{y};\kappa)
 g(\boldsymbol{y})\,\mathrm d\ell(\boldsymbol y)
\pm \frac12\boldsymbol\nu(\boldsymbol x)g(\boldsymbol x),
\quad \boldsymbol x\in S,
\]
and
\[
(\mathbf{curl}\,V)^\pm(\boldsymbol{x})
=
\int_S
\mathbf{curl}_{\boldsymbol{x}}G(\boldsymbol{x},\boldsymbol{y};\kappa)
 g(\boldsymbol{y})\,\mathrm d\ell(\boldsymbol y)
\mp \frac12\boldsymbol\tau(\boldsymbol x)g(\boldsymbol x).
\]
\end{proposition}

This follows from the free-space jump relations~\cite{ColtonKress1983} and $\mathbf{curl}\,V=(\nabla V)_\perp$.

To simplify the boundary integral equations, define
\begin{align}
\mathcal{S}_\kappa[g](\boldsymbol{x})
&=2\int_S G(\boldsymbol{x},\boldsymbol{y};\kappa)g(\boldsymbol{y})\,\mathrm d\ell(\boldsymbol y),\label{S}\\
\mathcal{K}_\kappa'[g](\boldsymbol{x})
&=2\int_S
\frac{\partial G(\boldsymbol{x},\boldsymbol{y};\kappa)}{\partial\boldsymbol\nu(\boldsymbol{x})}
 g(\boldsymbol{y})\,\mathrm d\ell(\boldsymbol y),\label{K'}\\
\mathcal{H}_\kappa[g](\boldsymbol{x})
&=2\int_S
\frac{\partial G(\boldsymbol{x},\boldsymbol{y};\kappa)}{\partial\boldsymbol\tau(\boldsymbol{x})}
 g(\boldsymbol{y})\,\mathrm d\ell(\boldsymbol y),\quad\boldsymbol x\in S.\label{H}
\end{align}

The fluid--solid transmission conditions also involve second-order normal traces of single-layer potentials. The above operators act componentwise on vector and matrix densities, and we write $D_\tau=\partial_{\boldsymbol\tau}$.
\begin{proposition}\label{order 2}
Let $0<\sigma<1$, let $\Gamma$ be a $C^{2,\sigma}$ periodic interface, and let $g\in C^{1,\sigma}_\alpha(S)$. Then
\begin{equation}\label{eq:bie:secondtraces}
\begin{aligned}
2\partial_{\boldsymbol\nu}(\nabla V)^\pm
 &=-\kappa^2\mathcal S_\kappa[\boldsymbol\nu\boldsymbol\nu^\top g]\boldsymbol\nu
   -\mathcal H_\kappa[D_\tau(g\boldsymbol\nu)]\\
 &\quad+\mathcal K'_\kappa[D_\tau(g\boldsymbol\tau)]
   \pm D_\tau(g\boldsymbol\tau),\\
2\partial_{\boldsymbol\nu}(\mathbf{curl}\,V)^\pm
 &=\kappa^2\mathcal S_\kappa[\boldsymbol\tau\boldsymbol\nu^\top g]\boldsymbol\nu
   +\mathcal H_\kappa[D_\tau(g\boldsymbol\tau)]\\
 &\quad+\mathcal K'_\kappa[D_\tau(g\boldsymbol\nu)]
   \pm D_\tau(g\boldsymbol\nu).
\end{aligned}
\end{equation}
\end{proposition}

\begin{proof}
We use tangential integration by parts for second-order single-layer traces~\cite{Dong2020}. Write
$\mathcal D_\kappa h=\int_S\partial_{\boldsymbol\nu(\boldsymbol y)}
G(\cdot,\boldsymbol y;\kappa)h(\boldsymbol y)\,\mathrm d\ell(\boldsymbol y)$.
Using $\nabla_{\boldsymbol x}G=-\nabla_{\boldsymbol y}G$,
$(\Delta+\kappa^2)G=0$, and tangential integration by parts gives, away from $\Gamma$,
\[
\begin{aligned}
\nabla V_\kappa g&=-\mathcal D_\kappa(g\boldsymbol\nu)
   +V_\kappa[D_\tau(g\boldsymbol\tau)],\\
\nabla\mathcal D_\kappa h&=\kappa^2V_\kappa(\boldsymbol\nu h)
   +\mathbf{curl}\,V_\kappa(D_\tau h).
\end{aligned}
\]
The second identity holds for scalar $h$ and is applied componentwise to vector densities. The kernel $G$ is $(-\alpha)$-quasi-periodic in the source variable, whereas the density is $\alpha$-quasi-periodic. Their product, and its products with the periodic vectors $\boldsymbol\nu,\boldsymbol\tau$, are therefore strictly periodic, so integration by parts produces no endpoint terms. Taking the normal derivative of the first identity and applying the second identity and Proposition~\ref{order 1} componentwise gives the first normal trace formula. Applying a clockwise rotation to this formula and using $\boldsymbol\nu_\perp=-\boldsymbol\tau$ and $\boldsymbol\tau_\perp=\boldsymbol\nu$ gives the second formula.
\end{proof}

\subsection{Coupled boundary integral system}\label{sec:bie:construction}
We first derive the integral system for general transmission data; the incident plane wave will then enter only through the right-hand side. Fix $0<\sigma<1$ and assume $\gamma\in C^{3,\sigma}_{\rm per}(\mathbb R)$. Let $\boldsymbol u$ and $w$ denote the elastic displacement on the solid side and the radiating acoustic field on the fluid side, respectively. Define the transmission operator
\begin{equation}\label{eq:bie:transmissionop}
\mathscr T_\gamma(\boldsymbol u,w)
=\begin{pmatrix}
\boldsymbol\nu\cdot\mathcal T\boldsymbol u+w\\
\boldsymbol\tau\cdot\mathcal T\boldsymbol u\\
\boldsymbol u\cdot\boldsymbol\nu-(\rho_a\omega^2)^{-1}\partial_{\boldsymbol\nu}w
\end{pmatrix}.
\end{equation}
Given $\boldsymbol f=(f_1,f_2,f_3)^\top\in(C^{2,\sigma}_\alpha(S))^3$, consider $\boldsymbol u,w$ satisfying the Navier and acoustic Helmholtz equations, the corresponding quasi-periodic and radiation conditions, and
\begin{equation}\label{eq:bie:generaltransmission}
\mathscr T_\gamma(\boldsymbol u,w)=\boldsymbol f.
\end{equation}
We call a weak solution pair of this problem a radiating local energy solution if both fields are locally in $H^1$ up to the interface on their respective sides; the transmission conditions are understood in the weak trace sense. The original scattering problem corresponds to $w=p^s$ and
\begin{equation}\label{eq:bie:incidentdata}
\boldsymbol f=\boldsymbol f^i
:=\left(-p^i,\ 0,\ \frac{\partial_{\boldsymbol\nu}p^i}{\rho_a\omega^2}\right)^\top.
\end{equation}

The two scalar potentials of the elastic displacement are determined by
\begin{equation}\label{eq:bie:recoverpot}
\phi=-\frac{1}{\kappa_p^2}\nabla\!\cdot\boldsymbol u,
\qquad
\psi=\frac{1}{\kappa_s^2}\operatorname{curl}\boldsymbol u.
\end{equation}
Taking the divergence and curl of the Navier equation and using vector identities gives
\begin{equation}\label{eq:bie:helmholtzdecomposition}
(\Delta+\kappa_p^2)\phi=0,\qquad
(\Delta+\kappa_s^2)\psi=0,\qquad
\boldsymbol u=\nabla\phi+\mathbf{curl}\,\psi.
\end{equation}
The elastic Rayleigh expansion shows that the potentials $\phi,\psi$ both satisfy the downward scalar radiation condition. Set $X_\sigma=(C^{1,\sigma}_\alpha(S))^3$ and $Y_\sigma=(C^{0,\sigma}_\alpha(S))^3$. For unknown densities $\boldsymbol g=(g_1,g_2,g_3)^\top\in X_\sigma$, define the candidate fields
\begin{align}
\phi(\boldsymbol x)&=(V_{\kappa_p}g_1)(\boldsymbol x),
\quad \boldsymbol x\in\Omega_0^-,\label{phi}\\
\psi(\boldsymbol x)&=(V_{\kappa_s}g_2)(\boldsymbol x),
\quad \boldsymbol x\in\Omega_0^-,\label{psi}\\
w(\boldsymbol x)&=(V_{\kappa_a}g_3)(\boldsymbol x),
\quad \boldsymbol x\in\Omega_0^+.\label{us}
\end{align}
These fields define the reconstruction map
\[
\mathcal R_\gamma\boldsymbol g
=\left(\nabla V_{\kappa_p}g_1+\mathbf{curl}\,V_{\kappa_s}g_2,\ V_{\kappa_a}g_3\right).
\]
The reconstructed fields automatically satisfy the governing equations, the quasi-periodic conditions, and the radiation conditions. Taking the lower traces of the elastic field and the upper traces of the acoustic field, and retaining the normalization in~\eqref{S}--\eqref{H}, we define
\begin{equation}\label{eq:bie:composition}
\mathcal B_\gamma=2\mathscr T_\gamma\mathcal R_\gamma.
\end{equation}
The candidate fields therefore satisfy the transmission conditions if and only if $\mathcal B_\gamma\boldsymbol g=2\boldsymbol f$. The trace formulas in Propositions~\ref{order 1} and~\ref{order 2} give the explicit system
\begin{align}
\label{BIE1}
2f_1 &= -2\mu \kappa_p^2 \boldsymbol{\nu}^\top \mathcal{S}_{\kappa_p}[\boldsymbol{\nu} \boldsymbol{\nu}^\top g_1] \boldsymbol{\nu} - 2\mu \boldsymbol{\nu}^\top \mathcal{H}_{\kappa_p}\left[ \frac{\partial (g_1 \boldsymbol{\nu})}{\partial \boldsymbol{\tau}} \right] + 2\mu \boldsymbol{\nu}^\top \mathcal{K}_{\kappa_p}'\left[ \frac{\partial (g_1 \boldsymbol{\tau})}{\partial \boldsymbol{\tau}} \right] \nonumber \\
&\quad + 2\mu \kappa_s^2 \boldsymbol{\nu}^\top \mathcal{S}_{\kappa_s}[\boldsymbol{\tau} \boldsymbol{\nu}^\top g_2] \boldsymbol{\nu} + 2\mu \boldsymbol{\nu}^\top \mathcal{H}_{\kappa_s}\left[ \frac{\partial (g_2 \boldsymbol{\tau})}{\partial \boldsymbol{\tau}} \right] + 2\mu \boldsymbol{\nu}^\top \mathcal{K}_{\kappa_s}'\left[ \frac{\partial (g_2 \boldsymbol{\nu})}{\partial \boldsymbol{\tau}} \right] \nonumber \\
&\quad - \lambda \kappa_p^2 \mathcal{S}_{\kappa_p}[g_1] + \mathcal{S}_{\kappa_a}[g_3] - 2\mu \boldsymbol{\nu} \cdot \frac{\partial (g_1 \boldsymbol{\tau})}{\partial \boldsymbol{\tau}} - 2\mu \boldsymbol{\nu} \cdot \frac{\partial (g_2 \boldsymbol{\nu})}{\partial \boldsymbol{\tau}}, \\
\label{BIE2}
2f_2 &= -2\mu\kappa_p^2 \boldsymbol{\tau}^\top \mathcal{S}_{\kappa_p}[\boldsymbol{\nu} \boldsymbol{\nu}^\top g_1] \boldsymbol{\nu} -2\mu \boldsymbol{\tau}^\top \mathcal{H}_{\kappa_p}\left[ \frac{\partial (g_1 \boldsymbol{\nu})}{\partial \boldsymbol{\tau}} \right] +2\mu \boldsymbol{\tau}^\top \mathcal{K}_{\kappa_p}'\left[ \frac{\partial (g_1 \boldsymbol{\tau})}{\partial \boldsymbol{\tau}} \right] \nonumber \\
&\quad + 2\mu\kappa_s^2 \boldsymbol{\tau}^\top \mathcal{S}_{\kappa_s}[\boldsymbol{\tau} \boldsymbol{\nu}^\top g_2] \boldsymbol{\nu} +2\mu \boldsymbol{\tau}^\top \mathcal{H}_{\kappa_s}\left[ \frac{\partial (g_2 \boldsymbol{\tau})}{\partial \boldsymbol{\tau}} \right] +2\mu \boldsymbol{\tau}^\top \mathcal{K}_{\kappa_s}'\left[ \frac{\partial (g_2 \boldsymbol{\nu})}{\partial \boldsymbol{\tau}} \right] \nonumber \\
&\quad -\mu \kappa_s^2 \mathcal{S}_{\kappa_s}[g_2] - 2\mu\boldsymbol{\tau} \cdot \frac{\partial (g_1 \boldsymbol{\tau})}{\partial \boldsymbol{\tau}} -2\mu \boldsymbol{\tau} \cdot \frac{\partial (g_2 \boldsymbol{\nu})}{\partial \boldsymbol{\tau}}, \\
\label{BIE3}
2f_3 &= \mathcal{K}_{\kappa_p}'[g_1] + \mathcal{H}_{\kappa_s}[g_2] - \frac{\mathcal{K}_{\kappa_a}'[g_3]}{\rho_a \omega^2} - g_1 - \frac{g_3}{\rho_a \omega^2}.
\end{align}

The system involves only these integral operators and local tangential differential terms, avoiding direct quadrature of hypersingular kernels.

The classical H\"older mapping properties of layer potentials give
\[
\mathcal S_\kappa:C^{0,\sigma}_\alpha(S)\to C^{1,\sigma}_\alpha(S),
\qquad
\mathcal K_\kappa',\mathcal H_\kappa:C^{0,\sigma}_\alpha(S)\to C^{0,\sigma}_\alpha(S).
\]
Together with the regularity of the geometric coefficients and tangential differential terms, these imply
\begin{equation}\label{eq:bie:Bop}
\mathcal B_\gamma\in\mathcal L(X_\sigma,Y_\sigma),
\qquad
\mathcal B_\gamma\boldsymbol g=2\boldsymbol f.
\end{equation}
Choosing the right-hand side in~\eqref{eq:bie:incidentdata} gives the integral system for the direct scattering problem.

\subsection{Equivalence and unique solvability}\label{sec:BIEequiv}
To recover the single-layer densities from a solution of the transmission problem, we first establish the following representation result for scalar radiating fields.

\begin{lemma}\label{lem:bie:densityreg}
Let $0<\sigma<1$ and $\gamma\in C^{2,\sigma}_{\rm per}(\mathbb R)$, and assume that $\kappa$ satisfies the non-Wood condition. Let $v$ be an $\alpha$-quasi-periodic Helmholtz solution in $\Omega^+$ or $\Omega^-$ that satisfies the corresponding Rayleigh condition and is locally in $H^1$ up to the interface on that side. If $v|_S\in C^{2,\sigma}_\alpha(S)$, then there is a unique $g\in H^{-1/2}_\alpha(S)$ such that $v=V_\kappa g$ in that domain. Moreover, $g\in C^{1,\sigma}_\alpha(S)$.
\end{lemma}

\begin{proof}
By Zhang and Zhang~\cite[Theorem~3.3]{ZhangZhang2018}, the single-layer trace operator
\begin{equation}\label{eq:bie:Siso}
\mathcal S_\kappa:H^{-1/2}_\alpha(S)\longrightarrow H^{1/2}_\alpha(S)
\end{equation}
is a bounded isomorphism. Set $h=v|_S$ and choose $g=2\mathcal S_\kappa^{-1}h$. Then $V_\kappa g$ and $v$ have the same boundary trace and satisfy the same quasi-periodic and radiation conditions on the chosen side. Uniqueness of the radiating Dirichlet problem~\cite[Lemma~3.1]{ZhangZhang2018} gives $v=V_\kappa g$, while injectivity of the single-layer trace operator gives uniqueness of the density.

To establish the classical regularity of the density, set $z=V_\kappa g$. Its local $H^1$ regularity on both sides, common Dirichlet trace $h$, and equation $\Delta z=-\kappa^2z\in L^2$ yield local $H^2$ regularity up to the interface~\cite[Theorem~8.12]{GilbargTrudinger2001}. The two-dimensional embedding then gives $z\in C^{0,\sigma}$, and local boundary Schauder regularity~\cite[Lemma~6.18]{GilbargTrudinger2001} yields $C^{2,\sigma}$ regularity on both sides. Hence
\[
g=(\partial_{\boldsymbol\nu}z)^+-(\partial_{\boldsymbol\nu}z)^-
\]
belongs to $C^{1,\sigma}$. A finite cover of one interface period, together with quasi-periodic compatibility, gives the desired conclusion.
\end{proof}

Applying this representation to the two displacement potentials and the acoustic field gives a one-to-one correspondence between solutions of the integral system and those of the transmission problem.

\begin{theorem}\label{thm:BIEequiv}
Let $0<\sigma<1$ and $\gamma\in C^{3,\sigma}_{\rm per}(\mathbb R)$, and suppose that condition~\eqref{ass:nonwood} holds. For any $\boldsymbol f\in(C^{2,\sigma}_\alpha(S))^3$, the reconstruction map $\mathcal R_\gamma$ gives a one-to-one correspondence between solutions in $X_\sigma$ of the integral system~\eqref{eq:bie:Bop} and radiating local energy solutions of the transmission problem~\eqref{eq:bie:generaltransmission}.
\end{theorem}

\begin{proof}
If $\boldsymbol g\in X_\sigma$ solves the integral system, the construction in the preceding subsection ensures that $\mathcal R_\gamma\boldsymbol g$ satisfies the governing equations, quasi-periodic conditions, and radiation conditions. Propositions~\ref{order 1} and~\ref{order 2} ensure that the required boundary traces exist, and~\eqref{eq:bie:composition} gives the transmission conditions~\eqref{eq:bie:generaltransmission}.

Conversely, let $(\boldsymbol u,w)$ be a local energy solution of the transmission problem. Write the interface conditions as
\[
\mathcal T\boldsymbol u=(f_1-w)\boldsymbol\nu+f_2\boldsymbol\tau,
\qquad
\partial_{\boldsymbol\nu}w=\rho_a\omega^2(\boldsymbol u\cdot\boldsymbol\nu-f_3).
\]
Starting from $H^1$, the coupling conditions and local boundary regularity for traction and Neumann problems yield $H^2$ and then $H^3$ regularity on both sides~\cite[Remark~6.3]{Yintao2016}. The two-dimensional embedding then gives $C^{1,\sigma}$ regularity. Substituting the fields into the interface conditions and applying local boundary Schauder regularity twice~\cite[Theorem~9.3]{AgmonDouglisNirenberg1964} raises the regularity of the boundary data and fields to $C^{2,\sigma}$ and $C^{3,\sigma}$, respectively. Thus, in a neighborhood $U$ of one interface period,
\[
\boldsymbol u\in C^{3,\sigma}(\overline{\Omega^-}\cap U)^2,
\qquad w\in C^{3,\sigma}(\overline{\Omega^+}\cap U).
\]
By~\eqref{eq:bie:recoverpot}, $\phi,\psi$ belong to $C^{2,\sigma}$ up to the interface on the solid side. Hence,
\begin{equation}\label{eq:bie:trace-reg}
\phi|_S,\quad\psi|_S,\quad w|_S\in C^{2,\sigma}_\alpha(S).
\end{equation}
Applying Lemma~\ref{lem:bie:densityreg} to the three scalar fields gives the densities $\boldsymbol g\in X_\sigma$ in~\eqref{phi}--\eqref{us}. Substitution into the transmission conditions yields $\mathcal B_\gamma\boldsymbol g=2\boldsymbol f$.

Finally, for any reconstructed potentials, taking the divergence and curl of the displacement recovers~\eqref{eq:bie:recoverpot}. Thus, the displacement uniquely determines $\phi,\psi$, and Lemma~\ref{lem:bie:densityreg} uniquely determines the single-layer densities of all three scalar fields. Hence $\mathcal R_\gamma$ is a bijection between the two solution sets.
\end{proof}

By this one-to-one correspondence, unique solvability of the transmission problem for given $\boldsymbol f$ implies that the integral system has a unique solution in $X_\sigma$. In particular, the frequencies chosen in Section~\ref{sec2} ensure a unique density solution for the original scattering problem. This result provides a common boundary integral framework for the original scattering problem and for the inhomogeneous derivative transmission problems considered later in the smooth setting.

\section{Nystr\"om discretization}\label{sec4}
To apply periodic quadrature to the integral system~\eqref{eq:bie:Bop}, define the pullback and phase-removal operator $\mathcal U_\alpha$ and the periodic tangential differential operator $\mathscr D_\alpha$ by
\begin{equation}\label{eq:quad:phase-transform}
(\mathcal U_\alpha f)(t)=\mathrm e^{-\mathrm i\alpha t}f(\boldsymbol\xi(t)),
\qquad
\mathscr D_\alpha=Q^{-1}(\partial_t+\mathrm i\alpha).
\end{equation}
Then $\widehat f=\mathcal U_\alpha f$ is periodic and $\mathcal U_\alpha\partial_{\boldsymbol\tau}\mathcal U_\alpha^{-1}=\mathscr D_\alpha$. For a periodic geometric coefficient $a$, we have
\[
\mathscr D_\alpha(a\widehat g)
=Q^{-1}\{a(\widehat g'+\mathrm i\alpha\widehat g)+a'\widehat g\}.
\]
Set
\begin{equation}\label{eq:quad:periodic-kernels}
\begin{aligned}
M(t,s;\kappa)&=2G(\boldsymbol\xi(t),\boldsymbol\xi(s);\kappa)
                         \mathrm e^{\mathrm i\alpha(s-t)},\\
K(t,s;\kappa)&=2(\gamma'(t),-1)^\top\cdot
 \nabla_{\boldsymbol x}G(\boldsymbol\xi(t),\boldsymbol\xi(s);\kappa)
                         \mathrm e^{\mathrm i\alpha(s-t)},\\
H(t,s;\kappa)&=2(1,\gamma'(t))^\top\cdot
 \nabla_{\boldsymbol x}G(\boldsymbol\xi(t),\boldsymbol\xi(s);\kappa)
                         \mathrm e^{\mathrm i\alpha(s-t)}.
\end{aligned}
\end{equation}
All three kernels are $2\pi$-periodic in both variables $t,s$. The periodic operators defined by
$\widehat{\mathcal A}=\mathcal U_\alpha\mathcal A\mathcal U_\alpha^{-1}$
are
\[
\begin{aligned}
\widehat{\mathcal S}_\kappa h(t)
 &=\int_0^{2\pi}M(t,s;\kappa)h(s)Q(s)\,\mathrm ds,\\
\widehat{\mathcal K}'_\kappa h(t)
 &=Q(t)^{-1}\int_0^{2\pi}K(t,s;\kappa)h(s)Q(s)\,\mathrm ds,\\
\widehat{\mathcal H}_\kappa h(t)
 &=Q(t)^{-1}\int_0^{2\pi}H(t,s;\kappa)h(s)Q(s)\,\mathrm ds.
\end{aligned}
\]
Principal values are defined as symmetric truncation limits after periodic extension. Removing the phase componentwise gives
\begin{equation}\label{BIEpa}
\widehat{\mathcal B}_\gamma\widehat{\boldsymbol g}=2\widehat{\boldsymbol f},
\qquad
\widehat{\mathcal B}_\gamma
=\mathcal U_\alpha\mathcal B_\gamma\mathcal U_\alpha^{-1}.
\end{equation}
For assembly, replace the operators in~\eqref{BIE1}--\eqref{BIE3} by their periodic counterparts while preserving the order of multiplication by the geometric coefficients.

Product quadrature is used for the logarithmic singularities, the Cauchy principal part is extracted using the periodic cotangent kernel, and the trapezoidal rule is applied to the smooth terms. To obtain smooth periodic remainders, we retain the true geometric distance and phase-removal factor in the local singular coefficients and localize them using a smooth cutoff~\cite{Kress1991,ZhangZhang2018}. Set
\[
\delta=\operatorname{wrap}_{[-\pi,\pi)}(t-s),\qquad
r_\delta=\sqrt{\delta^2+[\gamma(t)-\gamma(s)]^2},\qquad
E_\alpha(\delta)=\mathrm e^{-\mathrm i\alpha\delta}.
\]
Here $\delta$ is the representative modulo $2\pi$, and $r_\delta$ is the local distance after translation of adjacent periodic cells. Choose a smooth, periodic, even function $\chi$ that equals $1$ for $|\delta|\leq\pi/2$ and vanishes, together with all its derivatives, at $|\delta|=\pi$. In the numerical implementation, we use
\begin{equation*}
\chi(\delta)=
\begin{cases}
1,&|\delta|\le\pi/2,\\
\dfrac{\exp[-1/(1-u)]}{\exp(-1/u)+\exp[-1/(1-u)]},
   &\pi/2<|\delta|<\pi,\quad u=\dfrac{|\delta|-\pi/2}{\pi/2},\\
0,&|\delta|=\pi.
\end{cases}
\end{equation*}

Writing $\ell(t-s)=\ln(4\sin^2((t-s)/2))$, we decompose the kernels as
\begin{equation}\label{eq:quad:split}
\begin{aligned}
M(t,s;\kappa)&=M_1(t,s;\kappa)\ell(t-s)+M_2(t,s;\kappa),\\
K(t,s;\kappa)&=K_1(t,s;\kappa)\ell(t-s)+K_2(t,s;\kappa),\\
H(t,s;\kappa)&=\frac{1}{2\pi}\cot\frac{s-t}{2}
                  +H_1(t,s;\kappa)\ell(t-s)+H_2(t,s;\kappa),
\end{aligned}
\end{equation}
where the logarithmic coefficients are
\begin{equation}\label{eq:quad:coefficients}
\begin{aligned}
M_1(t,s;\kappa)&=-\frac{\chi(\delta)E_\alpha(\delta)}{2\pi}J_0(\kappa r_\delta),\\
K_1(t,s;\kappa)&=\frac{\kappa\chi(\delta)E_\alpha(\delta)}{2\pi}
 [\gamma(s)-\gamma(t)+\gamma'(t)\delta]\frac{J_1(\kappa r_\delta)}{r_\delta},\\
H_1(t,s;\kappa)&=\frac{\kappa\chi(\delta)E_\alpha(\delta)}{2\pi}
 [\delta+\gamma'(t)(\gamma(t)-\gamma(s))]\frac{J_1(\kappa r_\delta)}{r_\delta}.
\end{aligned}
\end{equation}
Here $J_0,J_1$ are the Bessel functions of the first kind of orders zero and one, respectively. The ratio $J_1(\kappa r_\delta)/r_\delta$ is defined on the diagonal by its limit $\kappa/2$, and the remainders $M_2,K_2,H_2$ are determined by~\eqref{eq:quad:split}.

For $\gamma\in C^\infty_{\rm per}$, all coefficients and remainders are smooth periodic functions. Indeed, $\chi=1$ near the diagonal, where~\eqref{eq:quad:coefficients} gives precisely the free-space logarithmic singular coefficients after phase removal. Both $J_0(\kappa r_\delta)$ and $J_1(\kappa r_\delta)/r_\delta$ are smooth functions of $r_\delta^2$, and the following ratio has a smooth positive extension,
\[
\frac{r_\delta^2}{4\sin^2(\delta/2)}
\longrightarrow 1+[\gamma'(t)]^2
\quad(s\to t).
\]
The remainders are therefore smooth after subtraction of the logarithmic terms and the cotangent principal part. Away from the diagonal, the original kernels are smooth, while the vanishing of $\chi$ to all orders at the periodic endpoints ensures smooth periodic extension.

Denote the local regular part and its spatial derivatives by
\[
G_\kappa^{\mathrm{reg}}(\boldsymbol x,\boldsymbol y)
=G(\boldsymbol x,\boldsymbol y;\kappa)-G_0(\boldsymbol x,\boldsymbol y;\kappa),\qquad
r_j(t;\kappa)=\left.\partial_{x_j}G_\kappa^{\mathrm{reg}}(\boldsymbol x,\boldsymbol y)\right|_{\boldsymbol x=\boldsymbol y=\boldsymbol\xi(t)}.
\]
Taking $s\to t$ gives the diagonal values
\[
\begin{aligned}
M_1(t, t; \kappa) &= -\frac{1}{2\pi}, \quad K_1(t, t; \kappa) = 0, \quad H_1(t, t; \kappa) = 0, \\
M_2(t, t; \kappa) &=  \frac{\mathrm{i}}{2} - \frac{\gamma_{e}}{\pi} - \frac{1}{\pi}\ln\left( \frac{\kappa\sqrt{1 + [\gamma'(t)]^2}}{2} \right)  +2G_\kappa^{\mathrm{reg}}(\boldsymbol\xi(t),\boldsymbol\xi(t)), \\
K_2(t, t; \kappa) &= -\frac{1}{2\pi}\frac{\gamma''(t)}{1 + [\gamma'(t)]^2} +2\bigl[\gamma'(t)r_1(t;\kappa)-r_2(t;\kappa)\bigr],\\
H_2(t, t; \kappa) &= -\frac{1}{2\pi}\frac{\gamma''(t)\gamma'(t)}{1 + [\gamma'(t)]^2} + \frac{\mathrm{i}\alpha}{\pi}+2\bigl[r_1(t;\kappa)+\gamma'(t)r_2(t;\kappa)\bigr].
\end{aligned}
\]
Here $\gamma_e$ is the Euler--Mascheroni constant. The term $\mathrm i\alpha/\pi$ in $H_2$ is the finite part generated by multiplying the Cauchy principal part by the phase factor used to obtain periodicity, and $r_j$ denotes a spatial derivative of the local regular part.

The kernel decomposition above yields a collocation discretization. Let $t_i=\pi i/N$ ($0\le i<2N$) be the collocation points.
Set $h_N=\pi/N$ and $Q_N=\operatorname{diag}(Q(t_i))$.
Let $\odot$ denote entrywise multiplication, and let
$M_{\ell,\kappa},K_{\ell,\kappa},H_{\ell,\kappa}$ ($\ell=1,2$)
denote the matrices obtained by evaluating the corresponding kernels at $(t_i,t_j)$. The Green function $G$ and its first derivatives are computed using Hettlich's fast integral representation~\cite{Hettlich2002}. The logarithmic and cotangent quadrature weights are
\[
\begin{aligned}
R_{ij}&=-\frac{2\pi}{N}\sum_{m=1}^{N-1}
 \frac{\cos[m(t_i-t_j)]}{m}
 -\frac{\pi}{N^2}\cos[N(t_i-t_j)],\\
C_{ij}&=\begin{cases}
 \dfrac{2\pi}{N}\cot\dfrac{t_j-t_i}{2},&j-i\text{ is odd},\\
 0,&j-i\text{ is even}.
\end{cases}
\end{aligned}
\]
The corresponding Nystr\"om matrices are
\begin{equation}\label{eq:quad:nystrom-matrices}
\begin{aligned}
S_{\kappa,N}&=(R\odot M_{1,\kappa}+h_NM_{2,\kappa})Q_N,\\
K'_{\kappa,N}&=Q_N^{-1}(R\odot K_{1,\kappa}+h_NK_{2,\kappa})Q_N,\\
H_{\kappa,N}&=Q_N^{-1}
 \left(\frac{C}{2\pi}+R\odot H_{1,\kappa}+h_NH_{2,\kappa}\right)Q_N.
\end{aligned}
\end{equation}
The diagonal kernel values are taken from the limits above. The differentiation matrix for periodic trigonometric interpolation is
\[
(D_N)_{ij}=
\begin{cases}
\dfrac{(-1)^{i-j}}{2}\cot\dfrac{(i-j)\pi}{2N},&i\ne j,\\
0,&i=j.
\end{cases}
\]
After phase removal, the tangential differential operator $\mathscr D_\alpha$ is represented by $Q_N^{-1}(D_N+\mathrm i\alpha I)$.
To discretize the product derivative in its expanded form, define, for a geometric coefficient $a$,
\begin{equation}\label{eq:quad:product-derivative}
T_{\alpha,a}
=Q_N^{-1}\{\operatorname{diag}(a)(D_N+\mathrm i\alpha I)
                         +\operatorname{diag}(a')\}.
\end{equation}
Apply this formula componentwise to $a=\boldsymbol\nu,\boldsymbol\tau$, first differentiating the continuous geometric coefficients and then evaluating them at the nodes. Finite-dimensional trigonometric differentiation does not satisfy an exact product rule, so this expression cannot be replaced by $Q_N^{-1}(D_N+\mathrm i\alpha I)\operatorname{diag}(a)$.

\begingroup\emergencystretch=2em
Integral terms containing product derivatives are discretized by differentiating first and then integrating. Thus $\mathcal H_\kappa[D_\tau(ga)]$ and $\mathcal K_\kappa'[D_\tau(ga)]$ are discretized as $H_{\kappa,N}T_{\alpha,a}$ and $K'_{\kappa,N}T_{\alpha,a}$, respectively. After assembly according to~\eqref{BIEpa} and solution for the densities, field values are computed from the layer potentials. The same kernel decomposition is used for the boundary DtN maps and the shape derivative equations.\par
\endgroup

\section{Fr\'echet shape derivative and its boundary integral computation}\label{sec:inverse}
This section links the shape linearization of near-field observations at finitely many measurement points to the integral solver developed above. We first establish local Fr\'echet differentiability of the observation map with a quadratic remainder estimate and then, under stronger smoothness assumptions, compute the observation Jacobian needed for reconstruction by solving the derivative transmission problem.

\subsection{Shape differentiability of the near-field observation map}\label{sec:inv:frechetsetup}
Following the fixed-domain variational approach to shape analysis of fluid--solid scattering \cite{BarucqDerivative2014,KangJiang2017,Barucq2017}, we establish remainder estimates for the solution and point observations on a single periodic cell with exact truncation conditions. Choose $H_-<\gamma^{\min}$ and $H_+>b$, and set
\[
\begin{aligned}
D_\gamma^+&=\{0<x_1<2\pi,\ \gamma(x_1)<x_2<H_+\},\\
D_\gamma^-&=\{0<x_1<2\pi,\ H_-<x_2<\gamma(x_1)\}.
\end{aligned}
\]
Impose exact acoustic and elastic DtN conditions on the upper and lower truncation lines, respectively, and denote the corresponding operators by $\mathscr T_a^+$ and $\mathscr T_e^-$. They are determined mode by mode from the upward acoustic and downward elastic Rayleigh expansions and map $H^{1/2}_\alpha$ continuously into $H^{-1/2}_\alpha$ \cite[Section~VII]{Yintao2016}. The outward normal to the lower truncation line is $-\boldsymbol e_2$. The acoustic operator has Fourier multiplier $\mathrm i\zeta_n$, and the denominator $d_n=\alpha_n^2+\beta_n\eta_n$ in the elastic modal matrix satisfies $\operatorname{Re}d_n\geq\kappa_p^2>0$.

Let $H^1_{\alpha,\#}$ denote the $H^1$ space whose traces on the left and right boundaries satisfy the $\alpha$-quasi-periodic condition, and set
\[
\mathscr X_\gamma=H^1_{\alpha,\#}(D_\gamma^-)^2\times H^1_{\alpha,\#}(D_\gamma^+),
\qquad Y=(\boldsymbol u,p^s).
\]
The duality pairing is linear in its first argument and conjugate linear in its second. After integration by parts, the contributions from the periodic side boundaries cancel, and the truncated problem takes the form
\begin{equation}\label{eq:inv:variational}
\mathscr A_\gamma Y=\ell_\gamma\quad\text{in }\mathscr X_\gamma^*.
\end{equation}
Choose a smooth periodic function $\vartheta$ that equals $1$ in a neighborhood of the interface and $0$ in neighborhoods of both truncation lines and the measurement line. For real-valued $\eta\in W^{1,\infty}_{\rm per}(0,2\pi)$, set
\begin{equation}\label{eq:inv:Theta}
\boldsymbol V_\eta(x)=\vartheta(x)\eta(x_1)\boldsymbol e_2,\qquad
\Theta_\eta(x)=x+\boldsymbol V_\eta(x).
\end{equation}
For sufficiently small $\|\eta\|_{W^{1,\infty}}$, $\Theta_\eta$ is a bi-Lipschitz homeomorphism from $D_\gamma^\pm$ onto $D_{\gamma+\eta}^\pm$ and commutes with translation by $2\pi\boldsymbol e_1$. The pullback therefore preserves the quasi-periodic structure, while the DtN operators and measurement points remain fixed. Writing
$\widetilde Y_\eta=(\boldsymbol u_{\gamma+\eta},p^s_{\gamma+\eta})\circ\Theta_\eta$, we obtain
\begin{equation}\label{eq:inv:pullbackop}
\mathscr A(\eta)\widetilde Y_\eta=\ell(\eta),\qquad
\mathscr A(0)=\mathscr A_\gamma,\quad \ell(0)=\ell_\gamma.
\end{equation}
Set $F_\eta=I+D\boldsymbol V_\eta$, $J_\eta=\det F_\eta$, $B_\eta=F_\eta^{-1}$, $C_\eta=\operatorname{cof}F_\eta$, and $M_\eta=J_\eta B_\eta B_\eta^\top$, and write
$E_\eta(\boldsymbol w)=\operatorname{sym}(D\boldsymbol w B_\eta)$ and $d_\eta(\boldsymbol w)=\operatorname{tr}(D\boldsymbol w B_\eta)$.
A change of variables and Nanson's formula $\boldsymbol\nu_\eta\,ds_\eta=C_\eta\boldsymbol\nu\,ds$ give the pullback form
\begin{equation}\label{eq:inv:pullbackA}
\begin{aligned}
&\langle\mathscr A(\eta)(\boldsymbol u,p^s),(\boldsymbol v,z)\rangle\\
&\quad=\int_{D_\gamma^-}\!J_\eta\bigl(2\mu E_\eta(\boldsymbol u):\overline{E_\eta(\boldsymbol v)}
+\lambda d_\eta(\boldsymbol u)\overline{d_\eta(\boldsymbol v)}
-\rho_e\omega^2\boldsymbol u\cdot\overline{\boldsymbol v}\bigr)\,dx\\
&\qquad+\int_{D_\gamma^+}\!\bigl(M_\eta\nabla p^s\cdot\nabla\overline z
-\kappa_a^2J_\eta p^s\overline z\bigr)\,dx
-\langle\mathscr T_e^-\boldsymbol u,\boldsymbol v\rangle_{H_-}
-\langle\mathscr T_a^+p^s,z\rangle_{H_+}\\
&\qquad-\int_S p^s(C_\eta\boldsymbol\nu)\cdot\overline{\boldsymbol v}\,ds
-\rho_a\omega^2\int_S(\boldsymbol u\cdot C_\eta\boldsymbol\nu)\overline z\,ds.
\end{aligned}
\end{equation}
\begin{equation}\label{eq:inv:pullbackell}
\begin{aligned}
\ell(\eta)(\boldsymbol v,z)
={}&\int_S(p^i\circ\Theta_\eta)(C_\eta\boldsymbol\nu)\cdot\overline{\boldsymbol v}\,ds\\
&-\int_S\bigl((\nabla p^i)\circ\Theta_\eta\bigr)\cdot(C_\eta\boldsymbol\nu)\overline z\,ds.
\end{aligned}
\end{equation}
Here $\langle\cdot,\cdot\rangle_{H_\pm}$ denotes the $H^{-1/2}_\alpha$--$H^{1/2}_\alpha$ antidual pairing on the truncation lines. The operator $\mathscr A_\gamma$ is Fredholm of index zero \cite{Yintao2016}. Thus, when the reference direct problem is uniquely solvable, it is a bounded isomorphism from $\mathscr X_\gamma$ onto $\mathscr X_\gamma^*$.

\begin{theorem}\label{thm:inv:solutionFrechet}\label{cor:inv:observationFrechet}
Suppose that $\gamma\in W^{1,\infty}_{\rm per}(\mathbb R)$ and that the reference direct problem \textbf{(\ref{prob:bvp})} is uniquely solvable. Then, in a $W^{1,\infty}_{\rm per}$ neighborhood of $0$, the perturbed problem is uniquely solvable, and the pullback solution map $\eta\mapsto\widetilde Y_\eta\in\mathscr X_\gamma$ is Fr\'echet differentiable at $0$. For every real-valued $q\in W^{1,\infty}_{\rm per}(0,2\pi)$, the material derivative $\dot Y[q]=(\dot{\boldsymbol u}[q],\dot p^s[q])$ is uniquely determined by
\begin{equation}\label{eq:inv:materialderivative}
\mathscr A_\gamma\dot Y[q]=\ell'(0)[q]-\mathscr A'(0)[q]Y_0,
\qquad Y_0=(\boldsymbol u,p^s),
\end{equation}
and satisfies the quadratic remainder estimate
\begin{equation}\label{eq:inv:solutionremainder}
\|\widetilde Y_\eta-Y_0-\dot Y[\eta]\|_{\mathscr X_\gamma}
\leq C\|\eta\|_{W^{1,\infty}}^2.
\end{equation}
The near-field observation map is well defined in this neighborhood and is Fr\'echet differentiable at $\gamma$, with
\begin{align}
\mathcal F'(\gamma)[q]&=\bigl(\dot p^s[q](\boldsymbol x^{(m)})\bigr)_{m=0}^{M-1},\label{eq:inv:Fprime}\\
\|\mathcal F(\gamma+\eta)-\mathcal F(\gamma)-\mathcal F'(\gamma)[\eta]\|_2
&\leq C\|\eta\|_{W^{1,\infty}}^2.\label{eq:inv:Frechetremainder}
\end{align}
\end{theorem}
\begin{proof}
Near $\eta=0$, the geometric coefficients $J_\eta,B_\eta,C_\eta,M_\eta$ depend smoothly on $D\boldsymbol V_\eta$, and the remainders of their first-order expansions in $L^\infty$ are $O(\|\eta\|_{W^{1,\infty}}^2)$. The interface coefficients are given explicitly in terms of $\eta,\eta'$ and satisfy the same estimate, as do the compositions with the smooth incident field. Equations \eqref{eq:inv:pullbackA}--\eqref{eq:inv:pullbackell}, the trace theorem, and the Cauchy--Schwarz inequality imply that $\mathscr A,\ell$ are continuously Fr\'echet differentiable. Writing $r_A(\eta)=\mathscr A(\eta)-\mathscr A_\gamma-\mathscr A'(0)[\eta]$ and $r_\ell(\eta)=\ell(\eta)-\ell(0)-\ell'(0)[\eta]$, we obtain
\begin{equation}\label{eq:inv:operatorremainder}
\|r_A(\eta)\|_{\mathcal L(\mathscr X_\gamma,\mathscr X_\gamma^*)}
+\|r_\ell(\eta)\|_{\mathscr X_\gamma^*}
\leq C\|\eta\|_{W^{1,\infty}}^2.
\end{equation}

Shrink the perturbation neighborhood so that $\|\mathscr A_\gamma^{-1}(\mathscr A(\eta)-\mathscr A_\gamma)\|<1/2$. The Neumann series gives a uniformly bounded inverse for $\mathscr A(\eta)$, while \eqref{eq:inv:materialderivative} defines a bounded linear map $q\mapsto\dot Y[q]$. For $R_\eta=\widetilde Y_\eta-Y_0-\dot Y[\eta]$, subtracting the reference and derivative equations from the perturbed equation gives
\[
\mathscr A(\eta)R_\eta=r_\ell(\eta)-r_A(\eta)Y_0
-\bigl(\mathscr A(\eta)-\mathscr A_\gamma\bigr)\dot Y[\eta].
\]
The right-hand side is $O(\|\eta\|_{W^{1,\infty}}^2)$, which yields \eqref{eq:inv:solutionremainder}.

Around each measurement point, choose fixed neighborhoods $B_m'\Subset B_m$ that lie in the fluid interior for all small perturbations and on which $\vartheta=0$. After quasi-periodic extension of the fields, local norms across a periodic cell boundary remain controlled by the energy norms on at most two periodic cells. In $B_m$, both $\dot p^s[\eta]$ and the pressure remainder $r_\eta=\widetilde p^s_\eta-p^s-\dot p^s[\eta]$ satisfy the homogeneous Helmholtz equation. Interior regularity and the two-dimensional embedding $H^2\hookrightarrow C^0$ give
$|r_\eta(\boldsymbol x^{(m)})|\leq C\|r_\eta\|_{H^1(B_m)}\leq C\|\eta\|_{W^{1,\infty}}^2$. Combining these estimates over the finitely many measurement points proves the observation result.
\end{proof}

On compact subsets shared by the reference and perturbed domains on either side of the interface, the shape derivatives at fixed spatial points are obtained from the material derivatives as
\begin{equation}\label{eq:inv:eulerfrommaterial}
\boldsymbol U[q]=\dot{\boldsymbol u}[q]-D\boldsymbol u\,\boldsymbol V_q,
\qquad P[q]=\dot p^s[q]-\nabla p^s\cdot\boldsymbol V_q.
\end{equation}
Since the incident field is independent of the shape, the shape derivative of the total acoustic pressure is also $P[q]$. Near the measurement line, $\boldsymbol V_q=0$, so $\dot p^s[q]$ in \eqref{eq:inv:Fprime} can be replaced by $P[q]$.

\subsection{Transmission characterization and boundary integral computation of the shape derivative}\label{sec:inv:transmission}
We use the local transmission conditions for the derivative established by Barucq et al.\ \cite[Eqs.~(3.7)--(3.8)]{BarucqMoussaoui2018}, together with the periodic differentiability result of the preceding subsection, to obtain the corresponding characterization. Let $0<\sigma<1$, $\gamma\in C^{3,\sigma}_{\rm per}$, and $q\in C^{2,\sigma}_{\rm per}$, and set
\begin{equation}\label{eq:inv:geometry}
\gamma_\varepsilon=\gamma+\varepsilon q,\qquad \boldsymbol h=q\boldsymbol e_2.
\end{equation}
Using the unit normal and tangent defined in Section~\ref{sec2}, we obtain
\begin{equation}\label{eq:inv:frame}
\delta_q\boldsymbol\nu=\frac{q'}{Q^2}\boldsymbol\tau,\qquad
\delta_q\boldsymbol\tau=-\frac{q'}{Q^2}\boldsymbol\nu.
\end{equation}
The reference fields $\boldsymbol u,p=p^i+p^s$ have $C^{2,\sigma}$ regularity up to the interface from their respective sides. Hence the second-order traces below are well defined in the classical sense, and $(\boldsymbol U,P)$ given by \eqref{eq:inv:eulerfrommaterial} belongs to $\mathscr X_\gamma$. Define the derivative potentials by
\begin{equation}\label{eq:inv:derivativepotentials}
\Phi[q]=-\kappa_p^{-2}\nabla\!\cdot\!\boldsymbol U[q],\qquad
\Psi[q]=\kappa_s^{-2}\operatorname{curl}\boldsymbol U[q].
\end{equation}

\begin{proposition}\label{prop:inv:derivative}
Under the preceding smoothness conditions and the unique solvability assumption of Theorem~\ref{thm:inv:solutionFrechet}, the shape derivatives satisfy
\begin{equation}\label{eq:inv:pdes}
\begin{aligned}
\Delta P+\kappa_a^2P&=0&&\text{in }\Omega^+,\\
\Delta\Phi+\kappa_p^2\Phi=\Delta\Psi+\kappa_s^2\Psi&=0&&\text{in }\Omega^-.
\end{aligned}
\end{equation}
They retain the quasi-periodicity parameter $\alpha$ and satisfy the corresponding upward or downward Rayleigh radiation conditions. On the reference interface,
\begin{equation}\label{eq:inv:vectorbc}
\mathcal T\boldsymbol U+P\boldsymbol\nu=-\boldsymbol R_1[q],\qquad
\rho_a\omega^2\boldsymbol U\cdot\boldsymbol\nu-\partial_{\boldsymbol\nu}P=-R_2[q],
\end{equation}
where
\begin{align}
\boldsymbol R_1[q]={}&\bigl((\boldsymbol h\cdot\nabla)\boldsymbol\sigma(\boldsymbol u)\bigr)\boldsymbol\nu
+\boldsymbol\sigma(\boldsymbol u)\delta_q\boldsymbol\nu
+(\boldsymbol h\cdot\nabla p)\boldsymbol\nu+p\,\delta_q\boldsymbol\nu,
\label{eq:inv:R1}\\
R_2[q]={}&\rho_a\omega^2\bigl\{(D\boldsymbol u\,\boldsymbol h)\cdot\boldsymbol\nu
+\boldsymbol u\cdot\delta_q\boldsymbol\nu\bigr\}\nonumber\\
&-\bigl\{\boldsymbol h^\top D^2p\,\boldsymbol\nu+\nabla p\cdot\delta_q\boldsymbol\nu\bigr\}.
\label{eq:inv:R2}
\end{align}
The right-hand sides are determined by the one-sided traces of the reference fields. This radiating transmission problem is uniquely solvable, and its acoustic pressure evaluated at the measurement points gives $\mathcal F'(\gamma)[q]$.
\end{proposition}
\begin{proof}
Taking difference quotients on fixed interior subdomains and using Theorem~\ref{thm:inv:solutionFrechet} and \eqref{eq:inv:eulerfrommaterial} yields the Helmholtz equation for $P$ and the Navier equation for $\boldsymbol U$. Taking the divergence and curl then gives \eqref{eq:inv:pdes} and $\boldsymbol U=\nabla\Phi+\mathbf{curl}\,\Psi$. The perturbation map preserves the quasi-periodic structure and equals the identity near the truncation lines, so the quasi-periodicity parameter and radiation conditions are unchanged.

The local transmission conditions can be differentiated in the weak flux-trace sense. Denote the displacement component of $\widetilde Y_{\varepsilon q}$ by $\widetilde{\boldsymbol u}_\varepsilon$, set $\widetilde p_\varepsilon=(p^i+p^s_{\gamma+\varepsilon q})\circ\Theta_{\varepsilon q}$, and define the pullback fluxes
$\boldsymbol S_\varepsilon=\bigl(2\mu E_{\varepsilon q}(\widetilde{\boldsymbol u}_\varepsilon)+\lambda d_{\varepsilon q}(\widetilde{\boldsymbol u}_\varepsilon)I\bigr)C_{\varepsilon q}$ and
$\boldsymbol g_\varepsilon=M_{\varepsilon q}\nabla\widetilde p_\varepsilon$.
By the Piola identity, the governing equations, and fixed-domain differentiability, both fluxes and their divergences are differentiable in $L^2$. Continuity of the $H(\operatorname{div})$ normal trace therefore allows differentiation in $H^{-1/2}_\alpha(S)$ of
\[
\boldsymbol S_\varepsilon\boldsymbol\nu+\widetilde p_\varepsilon C_{\varepsilon q}\boldsymbol\nu=0,
\qquad
\boldsymbol g_\varepsilon\cdot\boldsymbol\nu
=\rho_a\omega^2\widetilde{\boldsymbol u}_\varepsilon\cdot C_{\varepsilon q}\boldsymbol\nu.
\]
Substituting \eqref{eq:inv:eulerfrommaterial} and the variation of the unit normal, and using the original transmission conditions to eliminate the terms multiplied by the variation of the arclength element, gives \eqref{eq:inv:vectorbc}--\eqref{eq:inv:R2}. This agrees with the established local formulas cited above. Since the right-hand sides are $C^{0,\sigma}$ traces, local transmission regularity gives one-sided $C^{1,\sigma}$ regularity of $\boldsymbol U,P$, so the conditions also hold in the classical trace sense. The difference of any two solutions satisfies the reference homogeneous transmission problem and therefore vanishes by unique solvability of the reference problem.
\end{proof}

Write $\mathcal G_\gamma[q]=\mathcal F'(\gamma)[q]$. The normal component of the vertical velocity $q\boldsymbol e_2$ is $-q/Q$. Its tangential component produces only a periodic reparameterization and leaves observations at fixed points unchanged, in agreement with the Hadamard normal structure \cite{LaurainSturm2016}.

For boundary integral computation, assume further that $\gamma,q\in C^\infty_{\rm per}$. The regularity bootstrap in Section~\ref{sec3} gives one-sided regularity of every finite order for the reference fields. Set
\begin{equation}\label{eq:inv:B}
B_1=\boldsymbol\nu\cdot\boldsymbol R_1,\qquad
B_2=\boldsymbol\tau\cdot\boldsymbol R_1,\qquad B_3=R_2/(\rho_a\omega^2).
\end{equation}
Then the derivative transmission conditions become
\begin{equation}\label{eq:inv:scalarbc}
\mathscr T_\gamma(\boldsymbol U,P)=-(B_1,B_2,B_3)^\top.
\end{equation}
These right-hand sides are linear in $q,q'$. Their vector form involves only one-sided derivatives of $\boldsymbol u,p$ up to order two. After substituting $\boldsymbol u=\nabla\phi+\mathbf{curl}\,\psi$, numerical assembly requires Cartesian derivatives of $\phi,\psi$ up to order three and of $p$ up to order two.

Using the isomorphism property \eqref{eq:bie:Siso} of the single-layer trace operator and the first-order jump relations, define, for $\kappa\in\{\kappa_a,\kappa_p,\kappa_s\}$, the two-sided radiating DtN maps
\begin{equation}\label{eq:inv:dtn}
\Lambda_\kappa^\pm=(\mathcal K_\kappa'\pm I)\mathcal S_\kappa^{-1}:
H^{1/2}_\alpha(S)\to H^{-1/2}_\alpha(S).
\end{equation}
Here $+$ and $-$ correspond to the fluid and solid sides, respectively, and the inverse is applied by solving a linear system. For a physical trace $f=\mathrm e^{\mathrm i\alpha t}\widehat f$, we have
$\partial_{\boldsymbol\tau}f=\mathrm e^{\mathrm i\alpha t}\mathscr D_\alpha\widehat f$.
Set $\mathcal R_{\kappa,j}^\pm f=\tau_j\partial_{\boldsymbol\tau}f+\nu_j\Lambda_\kappa^\pm f$. If $w$ is a smooth radiating Helmholtz solution on the corresponding side and $f=w|_S$, then
\begin{equation}\label{eq:inv:recursion}
(\partial_{x_j}w)|_S=\mathcal R_{\kappa,j}^\pm f,\qquad
(\partial_{x_{j_r}}\cdots\partial_{x_{j_1}}w)|_S
=\mathcal R_{\kappa,j_r}^\pm\cdots\mathcal R_{\kappa,j_1}^\pm f.
\end{equation}
At each step, tangential differentiation acts on the entire field trace from the preceding step, including the geometric coefficients arising from the normal and tangent vectors. In the implementation, periodic derivatives are computed by FFT, while commutativity of mixed derivatives and
$(\partial_{11}+\partial_{22}+\kappa^2)\partial^\beta w=0$ ($|\beta|\leq1$) serve as consistency diagnostics. The lower-side recursion is used for $\phi,\psi$. For $p^s$, the upper-side outgoing recursion is used, and the known incident-field derivatives are then added to obtain the traces of the total acoustic pressure and its derivatives.

By Proposition~\ref{prop:inv:derivative} and the unique solvability conclusion following Theorem~\ref{thm:BIEequiv}, the derivative fields have a unique single-layer potential density $\widetilde{\boldsymbol g}\in X_\sigma$ satisfying
\begin{equation}\label{eq:inv:linearizedbie}
\mathcal B_\gamma\widetilde{\boldsymbol g}=-2(B_1,B_2,B_3)^\top.
\end{equation}
For a fixed interface and frequency, the direct problem and the linearized problems for all smooth directions share one discrete left-hand-side matrix. Once this matrix has been factorized, each direction requires only an updated right-hand side and back substitution. Evaluating the derivative single-layer potential at the measurement points gives the near-field response. Below, $\mathcal G_{N,\gamma}$ denotes its Nystr\"om discretization, and the responses to basis directions in a finite-dimensional Fourier space form the shape Jacobian.

\section{Multifrequency reconstruction using shape derivatives}\label{sec:newton}
We combine the shape derivative with a regularized backtracking Gauss--Newton iteration. Wavenumber continuation and progressive expansion of the Fourier space \cite{BaoLiWu2012,WangLinLv2026} allow the reconstruction to proceed from low frequencies and low-order modes.

Let the observation frequencies be $\omega_1<\cdots<\omega_{n_f}$, with corresponding Fourier orders
$0\leq J_1\leq\cdots\leq J_{n_f}=J_{\max}$ and observation vectors $\boldsymbol d_m^\delta\in\mathbb C^M$. At the $m$th stage, the interface is updated in $V_{J_m}=\operatorname{span}\{1,\cos t,\sin t,\ldots,\cos(J_m t),\sin(J_m t)\}$, where $V_0=\operatorname{span}\{1\}$. Omitting the frequency index, we represent the interface as
\begin{equation}\label{eq:inv:fourier}
\gamma_{\boldsymbol c}(t)=c_0+\sum_{j=1}^J\{a_j\cos(jt)+b_j\sin(jt)\},
\end{equation}
with coefficient vector $\boldsymbol c=(c_0,a_1,b_1,\ldots,a_J,b_J)^\top\in\mathbb R^{2J+1}$. At the first frequency, the iteration starts from a prescribed initial guess. At subsequent frequencies, the final coefficients from the preceding frequency are retained and newly activated modes are initialized to zero.

Each stage uses the corresponding wavenumbers $\kappa_{a,m},\kappa_{p,m},\kappa_{s,m}$ and quasi-periodicity parameter $\alpha^{(m)}=\kappa_{a,m}\cos\theta$. The current frequency, Fourier order, and discretization parameter are denoted by $\omega,J,N$.

Let $q_\ell$ denote the Fourier basis functions in the order given above. In a parameter neighborhood where the reference direct problem is uniquely solvable, Theorem~\ref{thm:inv:solutionFrechet} and the chain rule give the continuous parameter map and its derivative as
\begin{equation}\label{eq:inv:continuous-parameter-map}
\mathcal F_J(\boldsymbol c):=\mathcal F(\gamma_{\boldsymbol c}),\qquad
D\mathcal F_J(\boldsymbol c)\boldsymbol v
=\mathcal F'(\gamma_{\boldsymbol c})\left[\sum_{\ell=1}^{2J+1}v_\ell q_\ell\right].
\end{equation}
Let $\mathcal F_N(\gamma)\in\mathbb C^M$ denote the near-field observations computed using $2N$ boundary nodes. The admissible parameter set $\mathcal U_{N,J}\subset\mathbb R^{2J+1}$ requires $\max_{t\in[0,2\pi]}\gamma_{\boldsymbol c}(t)<b$ and invertibility of the discrete direct-problem matrix at the current frequency. Define the discrete parameter map
\begin{equation}\label{eq:inv:parameter-map}
F_{N,J}:\mathcal U_{N,J}\longrightarrow\mathbb C^M,\qquad
F_{N,J}(\boldsymbol c):=\mathcal F_N(\gamma_{\boldsymbol c}).
\end{equation}
All differentiation and linearization below take place in open neighborhoods contained in $\mathcal U_{N,J}$. For fixed $J$, we abbreviate the map as $F_N(\boldsymbol c)$.

The responses $\mathcal G_{N,\gamma}[q]$ of the discrete linearized transmission problem define the complex matrix at $\boldsymbol c_k$ by
\begin{equation}\label{eq:inv:jacobian-columns}
A_k\boldsymbol e_\ell=\mathcal G_{N,\gamma_{\boldsymbol c_k}}[q_\ell],
\qquad \ell=1,\ldots,2J+1,
\end{equation}
where $\boldsymbol e_\ell$ is the standard basis vector, $A_k\in\mathbb C^{M\times(2J+1)}$, and $q_\ell$ corresponds to the normal velocity $-q_\ell/Q$.

Since the coefficients are real, we write the residual and this matrix in real form as
\begin{equation}\label{eq:inv:real}
\boldsymbol r_k=\begin{pmatrix}\operatorname{Re}(\boldsymbol d^\delta-F_N(\boldsymbol c_k))\\
\operatorname{Im}(\boldsymbol d^\delta-F_N(\boldsymbol c_k))\end{pmatrix},\qquad
\mathbf J_k=\begin{pmatrix}\operatorname{Re}A_k\\\operatorname{Im}A_k\end{pmatrix}.
\end{equation}
Here the continuous problem is first linearized with respect to the shape, and the resulting transmission problem is then discretized. Thus, $A_k$ approximates $D\mathcal F_J(\boldsymbol c_k)$, and $\mathbf J_k\in\mathbb R^{2M\times(2J+1)}$ is its real form. Both matrices are referred to below as the shape Jacobian.

The local linear model is $m_k(\boldsymbol v)=F_{N,J}(\boldsymbol c_k)+A_k\boldsymbol v$. Given $\beta>0$, the increment is determined by Tikhonov regularization,
\begin{equation}\label{eq:inv:step}
\Delta\boldsymbol c_k
:=\mathop{\arg\min}_{\boldsymbol v\in\mathbb R^{2J+1}}
\left\{\|\mathbf J_k\boldsymbol v-\boldsymbol r_k\|_2^2
+\beta^2\|L_J\boldsymbol v\|_2^2\right\}.
\end{equation}
The weights for the Fourier increment are
\begin{equation}\label{eq:inv:weights}
L_J=\operatorname{diag}(w_0,w_1,w_1,\ldots,w_J,w_J),\qquad
w_j=(1+j^2)^{s/2},\quad j=1,\ldots,J,
\end{equation}
where $w_0>0$ and $s>0$, so higher modes are penalized more strongly. Equation~\eqref{eq:inv:step} is solved through the corresponding augmented least-squares system.

Define the data misfit as $\mathcal E(\boldsymbol c)=\tfrac12\|\boldsymbol d^\delta-F_N(\boldsymbol c)\|_2^2$. Backtracking tests $\boldsymbol v=\chi^h\Delta\boldsymbol c_k$ for $h=0,\ldots,T_{\max}-1$ and $0<\chi<1$, comparing the predicted and actual reductions
\begin{align}
\operatorname{pred}_k(\boldsymbol v)
&=\boldsymbol r_k^\top\mathbf J_k\boldsymbol v
-\tfrac12\|\mathbf J_k\boldsymbol v\|_2^2,\label{eq:inv:pred}\\
\operatorname{ared}_k(\boldsymbol v)
&=\mathcal E(\boldsymbol c_k)-\mathcal E(\boldsymbol c_k+\boldsymbol v),\qquad
\varrho_k=\frac{\operatorname{ared}_k(\boldsymbol v)}
{\operatorname{pred}_k(\boldsymbol v)}.\label{eq:inv:gain}
\end{align}
The actual reduction is computed by solving the candidate direct problem. The reduction ratio compares only the data misfit terms; increment regularization is used solely to determine the direction.

Given a minimum clearance $d_{\min}>0$ from the measurement line and a slope bound $s_{\max}>0$, the candidate interface must satisfy
\begin{equation}\label{eq:inv:admissible}
\max_l\gamma_{\boldsymbol c}(t_l)<b-d_{\min},\qquad
\max_l|\gamma'_{\boldsymbol c}(t_l)|<s_{\max}.
\end{equation}
These bounds are checked at the $2N$ inversion nodes. A step is accepted only if the candidate direct problem is solved successfully, the increment and residual are finite, both reductions are finite and positive, and the finite reduction ratio satisfies $\varrho_k\geq\eta_{\mathrm{acc}}$. Otherwise, only the step length is reduced; the direction is not recomputed and the Jacobian is not reassembled. Upon acceptance, set $\boldsymbol c_{k+1}=\boldsymbol c_k+\boldsymbol v$ and retain the direct solution for the next iteration. The index $k$ counts only accepted steps; rejected trials do not increase $k$.

Let $\delta_{\mathrm{stop},m}\geq0$ denote the stopping scale at the $m$th frequency. When it is positive, the residual criterion is
\begin{equation}\label{eq:inv:discrepancy}
\|\boldsymbol d_m^\delta-F_{N_m,J_m}^{(m)}(\boldsymbol c_{m,k})\|_2
\leq\tau_{\mathrm d}\delta_{\mathrm{stop},m},\qquad\tau_{\mathrm d}>1,
\end{equation}
where the superscript $(m)$ denotes the frequency $\omega_m$. For noiseless data, set $\delta_{\mathrm{stop},m}=0$ and use
\begin{equation}\label{eq:inv:noiselessstop}
\frac{\|\boldsymbol d_m^0-F_{N_m,J_m}^{(m)}(\boldsymbol c_{m,k})\|_2}
{\|\boldsymbol d_m^0\|_2}<\epsilon_{\mathrm r}.
\end{equation}

Before backtracking, the unscaled direction is also checked against the small-step criterion
\begin{equation}\label{eq:inv:stepstop}
\frac{\|\Delta\boldsymbol c_{m,k}\|_2}
{\max(\|\boldsymbol c_{m,k}\|_2,1)}<\epsilon_{\mathrm s}.
\end{equation}

The current frequency stage ends when the residual or small-step criterion is met, the number of accepted steps reaches $K_{\max}$, or all $T_{\max}$ backtracking trials are rejected. The corresponding reason is recorded, and continuation proceeds to the next frequency with the current coefficients and the new wavenumbers and grid. An inadmissible candidate geometry, a failed candidate direct solve, or an invalid candidate residual causes only the current trial to be rejected. If the initial direct solve at a stage, the Jacobian calculation, or the direction calculation fails or produces a nonfinite result, the entire trajectory terminates, returning the available records and the frequency at which the failure occurred. Normal completion returns the final-frequency interface, the stage statuses, and the total number of accepted steps $K=\sum_{m=1}^{n_f}K_m$, where $K_m$ counts the accepted steps at each frequency.

The algorithm parameters are listed in Table~\ref{tab:inv:algorithmparameters}. The frequencies, Fourier orders, grids, and observation configuration are given in Section~\ref{sec:inv:reconstructionexperiments}. Algorithm~\ref{alg:inverse} summarizes the procedure.

\begin{table}[!htbp]\centering\small
\caption{Parameters of the inverse algorithm and their values.}\label{tab:inv:algorithmparameters}
\begin{tabular*}{\textwidth}{@{\extracolsep{\fill}}lp{0.56\textwidth}l@{}}\toprule
Parameter & Role & Value\\\midrule
$\beta$ & Fixed regularization parameter for coefficient increments & $10^{-3}$\\
$w_0,\ s$ & Constant-mode weight and Fourier weight exponent & $0.2,\ 1$\\
$\chi$ & Backtracking step length reduction factor & $0.5$\\
$\eta_{\mathrm{acc}}$ & Candidate step acceptance threshold & $10^{-4}$\\
$\tau_{\mathrm d}$ & Safety factor for the noise threshold & $1.1$\\
$\epsilon_{\mathrm r}$ & Relative residual threshold for noiseless data & $10^{-5}$\\
$\epsilon_{\mathrm s}$ & Relative small-step threshold for the unscaled direction & $10^{-6}$\\
$d_{\min},\ s_{\max}$ & Minimum measurement-line clearance and nodal slope bound & $0.2,\ 2$\\
$K_{\max}$ & Maximum accepted steps per frequency & $18$\\
$T_{\max}$ & Maximum total backtracking trials per iteration & $12$\\
\bottomrule\end{tabular*}
\end{table}

\begin{algorithm}[!htbp]\small
\caption{Multifrequency interface reconstruction by regularized backtracking Gauss--Newton}\label{alg:inverse}
\begin{algorithmic}[1]
\Require Data, stopping scales, and $\{\omega_m,J_m,N_m\}$ at all frequencies, initial coefficients, and algorithm parameters.
\State Set $K=0$. If an initial direct solve, Jacobian calculation, or direction calculation fails or gives a nonfinite result, return the available records and the trajectory termination status.
\For{$m=1,\ldots,n_f$}
\State Set the initial coefficients, wavenumbers, $\alpha^{(m)}$, and $J_m,N_m$ as described above. Set $k=0$ and solve the initial direct problem.
\While{the current frequency stage has not ended}
\State Check the residual criterion \eqref{eq:inv:discrepancy} or \eqref{eq:inv:noiselessstop} and the accepted-step limit. End the stage if either is reached.
\State Assemble $\mathbf J_k$ from the responses in all Fourier directions and solve \eqref{eq:inv:step} for $\Delta\boldsymbol c_k$.
\State If \eqref{eq:inv:stepstop} holds, end the stage by the small-step criterion.
\State Test $\chi^h\Delta\boldsymbol c_k$ for $h=0,\ldots,T_{\max}-1$, accepting or rejecting each trial as described above. If all trials are rejected, end the stage with backtracking failure.
\State Update the coefficients and the direct state, and set $k\gets k+1$.
\EndWhile
\State Record $K_m=k$ and the stage status, and set $K\gets K+K_m$.
\EndFor
\State \Return The final-frequency interface, $K$, and the records and exit statuses at all frequencies.
\end{algorithmic}
\end{algorithm}

\FloatBarrier

\section{Numerical experiments}\label{sec5}
This section tests the discretization accuracy of the direct problem and examines the reconstruction of smooth periodic interfaces by multifrequency continuation.

\subsection{Accuracy and convergence for the direct scattering problem}\label{sec:num:validation}
The baseline material and incidence parameters are
\begin{equation}\label{eq:inv:physical}
\lambda=\mu=\rho_e=\rho_a=1,\qquad c=3,\qquad\theta=\pi/3.
\end{equation}
We use acoustic wavenumbers $\kappa_a=0.5,3,6.25$, corresponding to frequencies $\omega=1.5,9,18.75$. Throughout this section, $2N$ and $M$ denote the numbers of boundary nodes and measurement points, respectively.

For the planar interface $\gamma(t)=0$, the exact solution has the form
\[
\phi=A_p\mathrm e^{\mathrm i(\alpha x_1-\beta_0x_2)},\qquad
\psi=A_s\mathrm e^{\mathrm i(\alpha x_1-\eta_0x_2)},\qquad
p^s=A_a\mathrm e^{\mathrm i(\alpha x_1+\zeta_0x_2)},
\]
where the coefficients satisfy the linear system derived from the transmission conditions,
\[
\begin{pmatrix}
-2\mu\beta_0^2-\lambda\kappa_p^2&-2\mu\alpha\eta_0&1\\
-2\mu\alpha\beta_0&2\mu\eta_0^2-\mu\kappa_s^2&0\\
\beta_0&\alpha&\zeta_0/(\rho_a\omega^2)
\end{pmatrix}
\begin{pmatrix}A_p\\A_s\\A_a\end{pmatrix}
=\begin{pmatrix}-1\\0\\\zeta_0/(\rho_a\omega^2)\end{pmatrix}.
\]
The three interface traces are compared with the analytical solution using the relative maximum-norm errors
\begin{equation}\label{eq:num:flaterror}
e_N^{\mathrm{flat}}(v)=\frac{\|v_N-v_{\mathrm{ex}}\|_\infty}{\|v_{\mathrm{ex}}\|_\infty},
\qquad v\in\{\phi,\psi,p^s\},
\end{equation}
where $v_{\mathrm{ex}}$ is the analytical trace and the norm is evaluated at the boundary nodes. For the nonplanar interface
\begin{equation}\label{eq:num:target}
\gamma_*(t)=0.15\mathrm e^{\sin t}+0.1\mathrm e^{\cos(2t+\pi/2)},
\end{equation}
consistency between successive discretizations is assessed by the relative difference on their common coarse-grid nodes,
\begin{equation}\label{eq:num:selferror}
E_N^{\mathrm{self}}(v)=\frac{\|v_N-v_{2N}|_N\|_\infty}{\|v_{2N}|_N\|_\infty},
\end{equation}
where $v$ again denotes the interface trace of one of the three fields.

In Table~\ref{tab:direct:new}, the self-convergence differences at $N=128$ are referenced to the solution at $N=256$. At this level, the planar-interface errors in all three traces at all three frequencies are at most $3.85\times10^{-13}$, and the nonplanar self-convergence differences are at most $2.06\times10^{-12}$. Higher frequencies require finer grids. Once the fields are sufficiently resolved, the errors decrease rapidly under refinement, exhibiting the rapid convergence characteristic of smooth periodic quadrature until the numerical precision plateau is reached.

\begin{table}[!htbp]\centering\footnotesize\setlength{\tabcolsep}{2pt}
\caption{Errors relative to the analytical solution for the planar interface and self-convergence differences for the nonplanar interface at different frequencies.}\label{tab:direct:new}
\begin{tabular*}{\textwidth}{@{\extracolsep{\fill}}rcccccc@{}}\toprule
&\multicolumn{3}{c}{\shortstack{Planar-interface errors relative\\to the analytical solution\\$e_N^{\mathrm{flat}}$}}&\multicolumn{3}{c}{\shortstack{Nonplanar self-convergence differences\\$E_N^{\mathrm{self}}$}}\\
\cmidrule(lr){2-4}\cmidrule(lr){5-7}
$N$&$\phi$&$\psi$&$p^s$&$\phi$&$\psi$&$p^s$\\\midrule
\multicolumn{7}{c}{$\kappa_a=0.5$}\\
4&$1.506\times10^{-4}$&$1.730\times10^{-4}$&$3.332\times10^{-4}$&$3.171\times10^{-1}$&$7.646\times10^{-1}$&$5.446\times10^{-1}$\\
8&$4.165\times10^{-6}$&$5.110\times10^{-6}$&$9.466\times10^{-6}$&$4.355\times10^{-1}$&$9.779\times10^{-1}$&$4.206\times10^{-1}$\\
16&$2.257\times10^{-8}$&$2.456\times10^{-10}$&$2.273\times10^{-8}$&$4.768\times10^{-4}$&$1.913\times10^{-3}$&$3.564\times10^{-4}$\\
32&$9.127\times10^{-11}$&$2.200\times10^{-11}$&$6.829\times10^{-11}$&$7.805\times10^{-8}$&$2.278\times10^{-7}$&$5.452\times10^{-8}$\\
64&$3.410\times10^{-12}$&$3.118\times10^{-12}$&$6.410\times10^{-12}$&$1.625\times10^{-10}$&$4.660\times10^{-10}$&$9.377\times10^{-11}$\\
128&$1.011\times10^{-13}$&$3.845\times10^{-13}$&$4.525\times10^{-14}$&$4.375\times10^{-13}$&$1.027\times10^{-12}$&$5.995\times10^{-14}$\\\midrule
\multicolumn{7}{c}{$\kappa_a=3$}\\
4&$1.708\times10^{-1}$&$1.149\times10^{-1}$&$4.466\times10^{-1}$&$1.238$&$1.837$&$2.583$\\
8&$3.515\times10^{-2}$&$1.232\times10^{-2}$&$4.916\times10^{-2}$&$2.409\times10^{-1}$&$5.890\times10^{-1}$&$2.075\times10^{-1}$\\
16&$9.459\times10^{-6}$&$3.775\times10^{-6}$&$1.627\times10^{-5}$&$7.617\times10^{-2}$&$5.197\times10^{-1}$&$5.960\times10^{-2}$\\
32&$2.344\times10^{-8}$&$8.078\times10^{-9}$&$3.587\times10^{-8}$&$2.415\times10^{-6}$&$9.484\times10^{-6}$&$7.350\times10^{-7}$\\
64&$9.170\times10^{-12}$&$6.248\times10^{-12}$&$3.052\times10^{-11}$&$1.961\times10^{-10}$&$9.267\times10^{-10}$&$1.319\times10^{-10}$\\
128&$2.844\times10^{-14}$&$1.276\times10^{-13}$&$1.653\times10^{-14}$&$7.866\times10^{-14}$&$3.259\times10^{-13}$&$4.161\times10^{-14}$\\\midrule
\multicolumn{7}{c}{$\kappa_a=6.25$}\\
4&$1.535$&$3.434\times10^{-1}$&$2.199$&$1.732$&$1.251$&$1.605$\\
8&$9.651\times10^{-2}$&$1.064\times10^{-1}$&$3.920\times10^{-1}$&$2.040$&$3.411$&$4.242$\\
16&$4.585\times10^{-2}$&$1.593\times10^{-2}$&$6.340\times10^{-2}$&$2.802\times10^{-1}$&$2.469\times10^{-1}$&$1.518\times10^{-1}$\\
32&$2.249\times10^{-7}$&$4.501\times10^{-8}$&$1.669\times10^{-7}$&$5.772\times10^{-3}$&$3.269\times10^{-2}$&$3.533\times10^{-3}$\\
64&$8.785\times10^{-11}$&$2.196\times10^{-11}$&$1.013\times10^{-10}$&$1.739\times10^{-8}$&$7.516\times10^{-8}$&$4.925\times10^{-9}$\\
128&$5.957\times10^{-14}$&$3.501\times10^{-13}$&$4.424\times10^{-14}$&$4.222\times10^{-13}$&$2.055\times10^{-12}$&$2.415\times10^{-13}$\\
\bottomrule\end{tabular*}
\par\smallskip
\end{table}

\subsection{Interface reconstruction by multifrequency continuation}\label{sec:inv:reconstructionexperiments}

We consider a baseline interface, an interface with fine corrugations, and the baseline interface with different material parameters. Using the baseline parameters and incidence angle in \eqref{eq:inv:physical}, we choose wavenumbers satisfying \eqref{ass:nonwood} and increasing Fourier orders,
\[
(\kappa_a,J)=(0.5,0),(1,1),(2.25,2),(3,3),(4.25,4),(5,5),(6.25,6).
\]
We set $\omega=c\kappa_a$ and distribute $M=20$ measurement points uniformly over one period. Each stage uses only the complex scattered acoustic pressure data at the current frequency.

The baseline interface is the composite exponential interface from Section~\ref{sec:num:validation},
\[
\gamma_*(t)=0.15\mathrm e^{\sin t}+0.1\mathrm e^{\cos(2t+\pi/2)}.
\]
The second interface has a sixth-harmonic corrugation and is given by
\begin{equation}\label{eq:num:targetroot}
\gamma_{\mathrm r,*}(t)=0.2\sqrt{-2\sin t+2\cos(2t+\pi/2)+8}-0.05\cos(6t).
\end{equation}
The first two cases use the baseline material parameters. The different-material case uses the baseline interface with $\lambda=4,\mu=3,\rho_e=2$. Both sets of material parameters satisfy $\lambda+2\mu\leq\rho_e c^2$. All three cases use the initial guess $\gamma_0=0.8$, measurement height $b=1.2$, and the same frequency continuation schedule and iteration parameters.

The noiseless observations are $\boldsymbol d^0=\mathcal F_{N_{\mathrm d}}(\gamma_*)$, with $2N_{\mathrm d}$ and $2N_{\mathrm i}$ nodes used for data generation and inversion, respectively. The baseline interface and the different-material case use $(2N_{\mathrm d},2N_{\mathrm i})=(128,64)$. The corrugated interface uses the same grids for the first six frequencies and refined grids $(256,128)$ at the final frequency. The grids are selected after refinement checks of the near-field data and weighted Jacobian $\mathbf J L_J^{-1}$, and are fixed across repeated experiments. At each frequency, noisy data are generated by
\begin{equation}\label{eq:inv:noise}
\boldsymbol d^\delta=\boldsymbol d^0+\epsilon_{\mathrm n}\|\boldsymbol d^0\|_2\frac{\boldsymbol z}{\|\boldsymbol z\|_2},\qquad
z_j=X_j+\mathrm iY_j,\qquad\delta_{\rm noise}=\epsilon_{\mathrm n}\|\boldsymbol d^0\|_2,
\end{equation}
where $X_j,Y_j$ are independent standard normal variables before normalization. The stopping scale is the norm of the injected noise, $\delta_{\rm stop}=\delta_{\rm noise}$, giving the threshold $\tau_{\mathrm d}\delta_{\rm noise}$. For a fixed seed and frequency, all cases and noise levels share the same normalized noise vector; vectors at different frequencies are sampled independently.

The geometric error and relative observation residual at the current frequency are defined as
\begin{equation}\label{eq:inv:errors}
E_\Gamma(\boldsymbol c)=\frac{\|\gamma_{\boldsymbol c}-\gamma_*\|_{L^2(0,2\pi)}}{\|\gamma_*\|_{L^2(0,2\pi)}},\qquad
R(\boldsymbol c)=\frac{\|F_{N_{\mathrm i}}(\boldsymbol c)-\boldsymbol d^\delta\|_2}{\|\boldsymbol d^\delta\|_2}.
\end{equation}
The interface error is evaluated at $8192$ periodic nodes, and the algorithm parameters are listed in Table~\ref{tab:inv:algorithmparameters}. The total number of accepted steps over the seven frequencies is denoted by $K$.

The noiseless inversion of the baseline interface accepts $21$ steps in total and ends by the small-step criterion at the final frequency, $\kappa_a=6.25$ and $J=6$. The final interface error of $0.117\%$ is close to the best $L^2$ approximation error of $0.111\%$ in the Fourier space of order six.

Each case is reconstructed $50$ times at noise levels of $1\%,5\%,10\%$ using seeds $2026,\ldots,2075$, giving $450$ trajectories in total. Table~\ref{tab:inv:noise} reports the final-frequency interface errors, total accepted steps, and counts of reconstructions meeting the noise threshold at the final frequency.

\begin{table}[!htbp]\centering\small\setlength{\tabcolsep}{5pt}
\caption{Statistics of repeated reconstructions with noisy data for the three cases.}\label{tab:inv:noise}
\footnotesize\setlength{\tabcolsep}{2pt}
\renewcommand{\arraystretch}{0.95}
\begin{tabular*}{\textwidth}{@{\extracolsep{\fill}}cccrc@{}}\toprule
Case&Noise&\shortstack{$E_\Gamma$ (\%)\\mean$\pm$SD}&$\overline K$&\shortstack{Noise threshold met\\at final frequency\\(count/total)}\\\midrule
\multirow{3}{*}{\shortstack{Baseline\\interface}}&1\%&$0.276\pm0.0418$&20.5&50/50\\
&5\%&$1.26\pm0.235$&14.4&50/50\\
&10\%&$2.52\pm0.475$&9.2&50/50\\\midrule
\multirow{3}{*}{\shortstack{Corrugated\\interface}}&1\%&$0.126\pm0.0240$&26.5&50/50\\
&5\%&$0.673\pm0.188$&23.0&50/50\\
&10\%&$1.52\pm0.300$&19.3&50/50\\\midrule
\multirow{3}{*}{\shortstack{Different-material\\case}}&1\%&$0.254\pm0.0874$&45.3&3/50\\
&5\%&$0.871\pm0.349$&43.9&50/50\\
&10\%&$1.69\pm0.620$&39.5&50/50\\
\bottomrule\end{tabular*}
\par\smallskip
\begin{minipage}{0.96\textwidth}\footnotesize
Each group contains $50$ reconstructions. SD is the sample standard deviation; means and SDs of $E_\Gamma$ are given to three significant digits. $\overline K$ is the mean total number of accepted steps over the seven frequencies. The last column counts reconstructions meeting the noise threshold at the final frequency, out of $50$.
\end{minipage}
\end{table}
The mean interface error increases with noise in all three cases, ranging from $1.52\%$ to $2.52\%$ at $10\%$ noise. The maximum error over all noisy samples is $3.84\%$. The corrugated interface has a mean error of only $0.126\%$ at $1\%$ noise, showing that the final Fourier modes can recover the fine-scale corrugations.

The only failures to meet the noise threshold at the final frequency occur in the $1\%$ noise group of the different-material case. In that group, only $3/50$ meet the threshold, while the remaining $47$ end by the small-step criterion. The mean interface error is still $0.254\%$.

\begin{figure}[!htbp]\centering
\includegraphics[width=\textwidth]{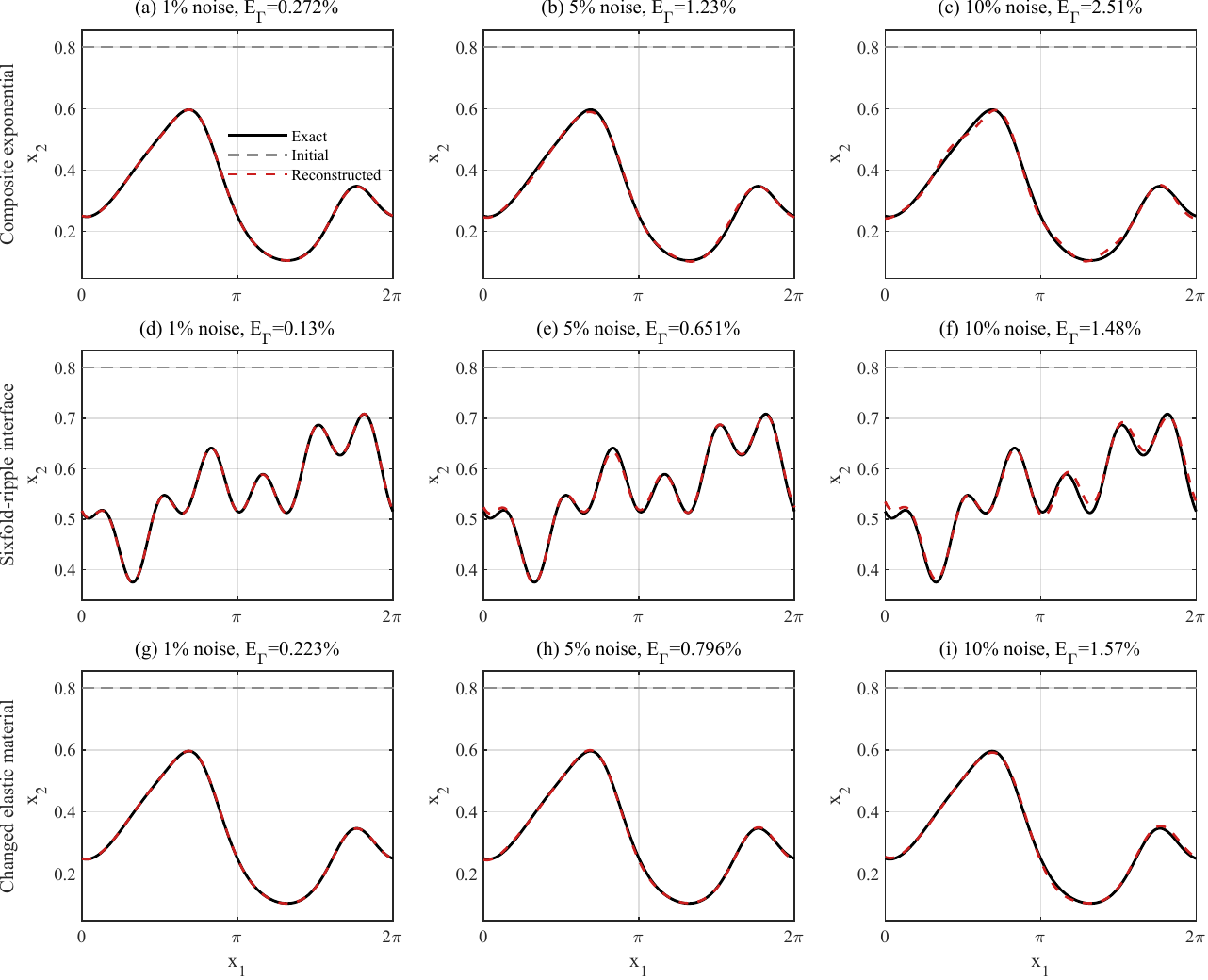}
\caption{Interface reconstructions from noisy data. The rows show the baseline interface (a--c), the corrugated interface (d--f), and the different-material case (g--i), respectively. The columns correspond to relative noise levels $\epsilon_{\mathrm n}$ of $1\%$, $5\%$, and $10\%$. Each panel shows the reconstruction with the $25$th smallest final-frequency relative interface error $E_\Gamma$ among the corresponding $50$ reconstructions.}\label{fig:inv:profiles}
\end{figure}

Figure~\ref{fig:inv:profiles} shows that both peaks and the trough are recovered for the baseline interface and the different-material case. The fine corrugations remain distinguishable at $10\%$ noise.

\FloatBarrier

\section{Conclusions and future work}\label{sec6}
We have developed an equivalent boundary integral formulation and a Nystr\"om discretization for periodic fluid--solid scattering. The shape derivative of the near-field observation map is characterized by a transmission problem, allowing the same framework to compute the Jacobian for multifrequency reconstruction. Numerical experiments support the direct solver's accuracy and show that the noiseless baseline reconstruction error is close to the best approximation error in the chosen Fourier space. Reconstructions from noisy data recover the main interface features across the tested geometries and material configurations. Future work may investigate reconstruction stability from finite near-field data and extend the method to three-dimensional doubly periodic fluid--solid interfaces.

\FloatBarrier

\ack
This work is supported by the National Key Laboratory of Computational Physics (Grant No. SYSM-2006-WDZC-13), the Science Challenge Project (Grant No. TZ2026004), and the National Natural Science Foundation of China (Grant Nos. 12271209 and 12671493).

\section*{Data availability statement}
The data cannot be made publicly available upon publication because they are not available in a format that is sufficiently accessible or reusable by other researchers. The data that support the findings of this study are available upon reasonable request from the authors.

\begingroup
\interlinepenalty=10000
\bibliographystyle{iopart-ip}
\bibliography{references}
\endgroup
\end{document}